\documentclass[10pt]{article}
\usepackage{graphicx}
\usepackage{amsmath,amssymb,amsthm,amsfonts}
\usepackage{amssymb}
\newtheorem{thm}{Theorem}[section]
\newtheorem{cor}[thm]{Corollary}
\newtheorem{lem}[thm]{Lemma}\newtheorem{prop}[thm]{Proposition}

\theoremstyle{definition}
\newtheorem{defn}{Definition}[section]
\theoremstyle{remark}
\newtheorem{rem}{Remark}[section]

\numberwithin{equation}{section}

\DeclareMathSymbol{\C}{\mathalpha}{AMSb}{"43}

\newcommand{\eps}{\varepsilon}

\newcommand{\lam}{\lambda}
\newcommand{\alp}{\alpha}
\newcommand{\sig}{\sigma}
\newcommand{\R}{{\mathbb{R}}}
\newcommand{\gam}{\gamma}

\def\endproof{\hfill$\blacksquare$\vspace{6pt}}

\newcommand{\bsub}{\begin{subequations}}
\newcommand{\esub}{\end{subequations}$\!$}
\begin{document}

\title{On the Partial Differential Equations of Electrostatic MEMS
Devices: Stationary Case}
\author{Nassif Ghoussoub\thanks{Partially supported by the Natural Science
and Engineering Research Council of Canada.}\ \ \ and \ Yujin
Guo\thanks{Partially supported by the Natural Science Foundation of
P. R. China (10171036) and by a U.B.C. Graduate Fellowship. }\\
Department of Mathematics, University of British Columbia,\\
Vancouver, B.C. Canada V6T 1Z2
\\ }
\date{}
\smallbreak \maketitle

\begin{abstract} We analyze the nonlinear elliptic problem  $\Delta u
=\frac{\lambda f(x)}{(1+u)^2}$ on a bounded domain $\Omega$ of
$\R^N$ with Dirichlet boundary conditions. This  equation models a
simple electrostatic Micro-Electromechanical System (MEMS) device
consisting of a thin dielectric elastic membrane with boundary
supported at $0$ above a rigid ground plate located at $-1$.   When
a voltage --represented here by $\lambda$-- is applied, the membrane
deflects towards the ground plate and a snap-through may occur when
it exceeds a certain critical value $\lambda^*$ (pull-in voltage).
This creates a so-called  ``pull-in instability" which greatly
affects the design of many devices. The mathematical model lends to
a nonlinear parabolic problem for the dynamic deflection of the
elastic membrane which will be considered in forthcoming papers
\cite{GG2} and \cite{GG3}.  For now, we focus on the stationary
equation where the challenge is to estimate $\lambda^*$ in terms of
material properties of the membrane, which can be fabricated with a
spatially varying dielectric permittivity profile $f$. Applying
analytical and numerical techniques, the existence of $\lambda^*$ is
established together with rigorous bounds. We show the existence of
at least one steady-state when $\lambda < \lambda^*$ (and when
$\lambda=\lambda^*$ in dimension $N< 8$) while none is possible for
$\lambda>\lambda^*$. More refined properties of steady states --such
as regularity, stability, uniqueness, multiplicity, energy estimates
and comparison results-- are shown to depend on the dimension of the
ambient space and on the  permittivity profile.

\end{abstract}

\vskip 0.2truein


\vskip 0.2truein

\tableofcontents

\section{Introduction}

Micro-Electromechanical Systems (MEMS) are often used to combine
electronics with micro-size mechanical devices in the design of
various types of microscopic machinery. MEMS devices have therefore
become key components of many commercial systems, including
accelerometers for airbag deployment in automobiles, ink jet printer
heads, optical switches and chemical sensors and so on (see for
example \cite{PB}). The simplicity and importance of this technique
have inspired numerous researchers to study mathematical models of
electrostatic-elastic interactions. The mathematical analysis of
these systems started in  the late 1960s with the pioneering work of
H.~C. Nathanson and his coworkers \cite{NN} who constructed and
analyzed a mass-spring model of electrostatic actuation, and offered
the first theoretical explanation of pull-in instability. At roughly
the same time, G.~I. Taylor \cite {T} studied the electrostatic
deflection of two oppositely charged soap films, and he predicted
that when the applied voltage was increased beyond a certain
critical voltage, the two soap films would touch together. Since
Nathanson and Taylor's seminal work, numerous investigators have
analyzed and developed mathematical models of electrostatic
actuation in attempts to understand further and control pull-in
instability. An overview of the physical phenomena of the
mathematical models associated with the rapidly developing field of
MEMS technology is given in \cite{PB}.

\begin{figure}[htbp]
\begin{center}
\includegraphics[width = 12cm, height = 5cm, clip]{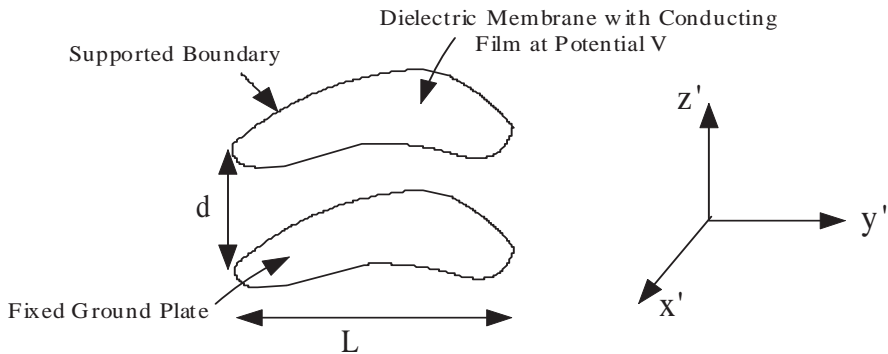}
\caption{{\em The simple electrostatic MEMS device.}}
\label{fig:fig1}
\end{center}
\end{figure}

The key component of many modern MEMS is the simple idealized
electrostatic device shown in Fig.~\ref{fig:fig1}. The upper part of
this device consists of a thin and deformable elastic membrane that
is held fixed along its boundary and which lies above a rigid
grounded plate. This elastic membrane is modeled as a dielectric
with a small but finite thickness. The upper surface of the membrane
is coated with a negligibly thin metallic conducting film. When a
voltage $V$ is applied to the conducting film, the thin dielectric
membrane deflects towards the bottom plate, and when $V$ is
increased beyond a certain critical value $V^*$ --known as pull-in
voltage--  the steady-state of the elastic membrane is lost, and
proceeds to touchdown or snap through at a finite time creating the
so-called pull-in instability.

A mathematical model of the physical phenomena, leading to a partial
differential equation for the dimensionless dynamic deflection of
the membrane, was derived and analyzed in \cite{FMP} and \cite{GPW}.
In the damping-dominated limit, and using a narrow-gap asymptotic
analysis, the dimensionless dynamic deflection $u=u(x,t)$ of the
membrane on a bounded domain $\Omega$ in $\mathbb{R}^2$, is found to
satisfy the following parabolic problem  \bsub \label{i:u}
\arraycolsep=1.5pt
\begin{eqnarray}
\frac{\partial u}{\partial t}& =& \Delta u - \frac{\lambda
f(x)}{(1+u)^2} \quad \quad {\rm for}
\quad x \in \Omega \,,\label{i:u1}\\
u (x,t) &=&0  \quad \quad \quad \quad \quad \quad \quad  {\rm for}
\quad x \in \partial\Omega \,,
\label{i:u2}\\
u(x,0) &=&0  \quad \quad \quad \quad \quad \quad \quad  {\rm for}
\quad x \in \Omega \,. \label{i:u3}
\end{eqnarray}
\esub An outline of the derivation of (\ref{i:u}) was given in
Appendix A of \cite{GPW}.  This initial condition in (\ref{i:u3})
assumes that the membrane is initially undeflected and the voltage
is suddenly applied to the upper surface of the membrane at time
$t=0$. The parameter $\lambda >0$ in (\ref{i:u1}) characterizes the
relative strength of the electrostatic and mechanical forces in the
system, and is given in terms of the applied voltage $V$ by
\begin{equation}
\lambda =\frac{\eps_0 V^2 L^2}{2T_e d^3} \,,\label{i:pullin}
\end{equation}
where $d$ is the undeflected gap size (see Fig.~\ref{fig:fig1}), $L$
is the length scale of the membrane, $T_e$ is the tension of the
membrane, and $\eps_0$ is the permittivity of free space in the gap
between the membrane and the bottom plate. In view of relation
$(\ref{i:pullin})$, we shall  use from now on the parameter $\lam $
and $\lam ^*$ to represent the applied voltage $V$ and pull-in
voltage $V^*$, respectively. Referred to as the {\em permittivity
profile}, $f(x)$ in $(1.1a)$ is defined by the ratio
\begin{equation}
     f(x) = \frac{\eps_0}{\eps_{2}(x)}  \,, \label{i:perm}
\end{equation}
where $\eps_{2}(x)$ is the dielectric permittivity of the thin
membrane.

There are several issues that must be considered in the actual
design of MEMS devices. Typically one of the primary goals is to
achieve the maximum possible stable deflection before touchdown
occurs, which is referred to as {\em pull-in distance} (cf.
\cite{GPW} and \cite{P1}). Another consideration is to increase the
stable operating range of the device by improving the pull-in
voltage $\lam ^*$ subject to the constraint that the range of the
applied voltage is limited by the available power supply. Such
improvements in the stable operating range is important for the
design of certain MEMS devices such as microresonators. One way of
achieving larger values of $\lam ^*$, while simultaneously
increasing the pull-in distance, is to use a voltage control scheme
imposed by an external circuit in which the device is placed
(cf.~\cite{PT}). This approach leads to a nonlocal problem for the
dynamic deflection of the membrane. A different approach studied in
\cite{P1} and \cite{GPW} is to introduce a spatially varying
dielectric permittivity $\eps_{2}(x)$ of the membrane. The idea is
to locate the region where the membrane deflection would normally be
largest under a spatially uniform permittivity, and then make sure
that a new dielectric permittivity $\eps_{2}(x)$ is largest --and
consequently  the profile $f(x)$ smallest-- in that region.

This latter approach requires the membrane having varying dielectric
properties, a framework investigated recently in \cite{P1} and
\cite{GPW}.  In \cite{P1} J.~Pelesko studied the steady-states of
$(1.1)$, when $f(x)$ is assumed to be bounded away from zero, i.e.,
\begin{equation}
   0 < C_0 \le f(x) \leq 1 \quad {\rm for\, all}\quad  x \in \Omega.
\label{i:fcon1}
\end{equation}
He established in this case an upper bound $\bar{\lambda}_1$ for
$\lambda^*$, and derived numerical results for the power-law
permittivity profile, from which the larger pull-in voltage and
thereby the larger pull-in distance, the existence and multiplicity
of the steady-states were observed. Recently, Y.~Guo, Z.~Pan and
M.~Ward studied in \cite{GPW} the dynamic behavior of $(1.1)$, which
is also of great practical interest. They considered a more general
class of profiles $f(x)$, where the membrane is allowed to be
perfectly conducting, i.e.,
\begin{equation}
0 \leq f(x) \leq 1  \quad {\rm for\, all}\quad  x \in \Omega \,,
\label{i:f}
\end{equation}
with $f(x)>0$ on a subset of positive measure of $\Omega$. By using
both analytical and numerical techniques, they obtained larger
pull-in voltage $\lam ^*$ and larger pull-in distance for different
classes of varying permittivity profiles. They obtained new upper
bounds $\bar{\lambda}_2$ for the pull-in voltage $\lam ^*$ that
correspond to the more general profiles $f(x)$ satisfying
$(\ref{i:f})$. Moreover, having estimated $\lam ^*$ numerically as
some saddle-node bifurcation value, they showed that $\lam ^*$ is
generally strictly smaller than both $\bar{\lambda}_1$ and
$\bar{\lambda}_2$ (cf. Table 1 of \cite{GPW}).

In this paper, we shall focus on the stationary deflection of the
elastic membrane, leaving the dynamic case to our forthcoming papers
\cite{GG2} and \cite{GG3}. For convenience, we shall set $v=-u$ in
such a way that our discussion will center on the following elliptic
problem
$$
\arraycolsep=1.5pt\begin{array}{lll}\arraycolsep=1.5pt
   \quad   \quad  \quad  \quad  \quad   -\Delta v &=& \displaystyle\frac{\lambda
f(x)}{(1-v)^2}    \,\,\quad x \in \Omega \,;\quad \\
   \hfill 0<&v&<1  \quad\quad   \quad   x \in \Omega \,;\\
 \hfill   v   &=& 0  \,\,\,\quad\quad\  \quad x \in
\partial\Omega.\quad
\end{array}\eqno{(S)_{\lam}}$$
We shall continue the investigation of optimal upper and lower
bounds for the pull-in voltage, and how they relate to the
permittivity profile $f$ which will be assumed to satisfy
(\ref{i:f}) throughout the paper unless mentioned otherwise. We
shall also discuss the issues of existence, multiplicity, and other
related properties of steady-states for $(S)_{\lam}$, and their
remarkable dependence on space-dimension.

This paper is organized as follows: In \S 2 we mainly show the
existence of a specific pull-in voltage and establish lower and
upper bound estimates. For that, we shall write  $\omega _{_N}$ for
the volume of the unit ball $B_1(0)$ in $\mathbb{R}^N$, and for any
bounded domain $\Gamma$ in $\R^n$, we denote by $\mu_{_\Gamma} $ the
first eigenvalue of $-\Delta$ on $H^1_0(\Gamma)$and by
$\phi_{_\Gamma}$ (resp., $\psi_{_\Gamma}$) the corresponding
positive eigenfunction normalized with $\int _{\Gamma
}\phi_{_\Gamma}dx=1$ (resp., $\sup_{x\in \Gamma}\psi_{_\Gamma}=1$).
We shall also associate to any domain $\Omega$ in $\R^N$  the
following parameter:
\begin{equation}
\nu_{_{\Omega}}=\sup\{\mu_{_\Gamma}H(\inf_\Omega \psi_\Gamma); \,
\hbox{\rm $\Gamma$ domain of $\R^N$, $\Gamma \supset \bar \Omega$}\}
\end{equation}
where $H$ is the function $H(t)=\frac{t(t+1+2\sqrt t)}{(t+1+\sqrt
t)^3}$.

We then prove the following lower estimates, the upper ones having
been established in \cite{P1} and \cite{GPW}.

\begin{thm}
  There exists a finite pull-in voltage $\lam ^*>0$ such that
\begin{enumerate}
\item If $0\leq \lam <\lam ^*$, there exists at least one solution for
$(S)_{ \lam}$;
\item If  $\lam >\lam ^*$, there is no solution for $(S)_{\lam}$.
\end{enumerate}
Moreover, we have the bounds
\[
\max \{\frac{\nu _{_\Omega}  }{\sup\limits_{x\in \Omega} f(x)},
\displaystyle\frac{8N}{27\sup\limits_{x\in \Omega}
f(x)}\Big(\displaystyle\frac{\omega _{_N}
}{|\Omega|}\Big)^{\frac{2}{N}}\}=:\underline\lam \le\lam ^*\leq
{\bar{\lambda}} :=\min\Big \{\displaystyle\frac{4\mu
_{_\Omega}}{27\inf\limits_{x\in \Omega }f(x)} ,
\displaystyle\frac{\mu_{_\Omega}}{3\int_{\Omega} f \phi_{_{\Omega}}
\, dx}\Big \}
       \]
      Furthermore, if $f(x)\equiv |x|^{\alp}$ on $\Omega$ with $\alp \ge
0$, then we have the more refined lower bound
\begin{equation}
  \lam _c(\alpha):=\frac{4(2+\alp)(N+\alp )}{27}
\Big(\displaystyle\frac{\omega
_{_N}}{|\Omega|}\Big)^{\frac{2+\alp}{N}}\le\lam ^*\,. \label{1:61}
\end{equation}

\end{thm}

In \S 2.3 we give some  numerical estimates on $\lam^*$ to compare
them with analytic bounds given in Theorem 1.1. Note that the upper
bound $\bar{\lambda}_1=\frac{4\mu _{_\Omega} }{27} \Big(\inf_{x\in
\Omega }f(x)\Big)^{-1}$ is relevant only when $f$ is bounded away
from $0$, while the upper bound $\bar{\lambda}_2=
\frac{\mu_{_\Omega}}{3} \Big( \int_{\Omega} f \phi_{_{\Omega}} \, dx
\Big)^{-1} $ is valid for all permittivity profiles.  In the case of
a uniform permittivity profile $f\equiv 1$ on $\Omega $, where
$\Omega $ is a strictly star-shaped domain containing $0$, we give a
more explicit upper bound $\bar{\lambda}_3$ of $\lam ^*$ in
Proposition 2.4.  In particular, we show that
$\bar{\lambda}_3=\frac{(N+2)^2}{8}$ is an upper bound in the case
where the domain is the unit ball $\Omega =B_1(0)
\subset\mathbb{R}^N$.

The issues of uniqueness and multiplicity of solutions for $(S)_{
\lam }$ with $0<\lam <\lam ^*$, and even mere existence for $(S)_{
\lam ^*}$ seem to be quite interesting. We address these problems
beginning in section \S 3
   by first considering minimal (positive) solutions of $(S)_{\lam}$
defined as follows.

\begin{defn} A solution $0<u_{\lam}(x)<1$ is said to be a minimal
(positive) solution of $(S)_{ \lam}$, if for any solution $0<u(x)<1$
of $(S)_{\lam}$ we have $u_{\lam}(x)\le u(x)$ in $\Omega $.
\end{defn}

Our main results in this direction can be stated as follows.

\begin{thm}Under the above assumptions, and with $\lam ^*$ as defined in
Theorem 1.1, there exists for any $\lam <\lam ^*$, a unique minimal
positive classical solution $u_{\lam}(x)$ of $(S)_{ \lam}$. It is
obtained as the limit of  the sequence $\{u_n(\lam ; x)\}$
constructed recursively as follows:
     $u_0\equiv 0$ in $\Omega $ and for each $n\geq 1$,
\begin{equation}\arraycolsep=1.5pt \begin {array}{lll}
-\Delta u_n &= \displaystyle\frac{\lambda f(x)}{(1-u_{n-1})^{2}} \,,
\quad &x \in \Omega \,;\\[3mm]
\quad 0\le &u_n<1 \,, \quad \quad\quad &x \in \Omega \,;
\\[3mm]
\quad \  u_n &= 0 \,,\quad \quad \quad \quad &x \in
\partial\Omega.
\end{array} \label{1:max}\end{equation}
Moreover, minimal solutions satisfy the following properties:
\begin{enumerate}
\item For each $x\in \Omega $, the function $\lambda \to u_{\lam}(x)$ is
strictly increasing and differentiable on $(0, \lambda^*)$;
\item If $1\le N<8$, then  there exists a constant
$0<C(N)<1$ such that $\parallel  u_{_\lam}
\parallel _{_{C(\Omega)}}\le C(N)$ for all
$\lam <\lam ^*$.
\end{enumerate}
\end{thm}

We refer to Lemma 3.6 in \S3.2 for a more general version of Theorem
1.2(2). Based on the results of Theorem 1.2, the existence and
related properties of minimal solutions at critical voltage $\lam
=\lam ^*$ will be studied in \S 3.3. More precisely, we shall
establish the following.

\begin{thm} If $1\le N<8$ then
$u_{_{\lam^*}}=\lim _{_{\lam \nearrow \lam ^*}}u_{_\lam}$ exists  in
the topology of $C^{2,\alp}(\bar {\Omega})$ with $0<\alp <1$, and
$u_{_{\lam ^*}}$ is the unique classical   solution of $(S)_{\lam
^*}$.
\end{thm}

\S 4 is devoted to the uniqueness and multiplicity of solutions
which remarkably depend again on the space-dimension.

\begin{thm} Under the above assumptions, with $\lam^*$ defined as in
Theorem 1.1,  we have:
\begin{enumerate}
\item If $N>2$, then for any $M>0$ there exists a voltage $0<\lam^*
_1(M)<\lam ^*$ such that for every $\lam \in (0,\lam^* _1(M))$,
there exists a unique positive solution for $(S)_\lambda$ --namely
the minimal solution $u_\lambda$-- that satisfies $\int _{\Omega}
|\frac{f}{(1-u)^3} |^{\frac{N}{2}}dx \leq M$;
\item If $1\le N<8$ then there exists $0<\lam ^*_2<\lam ^*$ such that
$(S)_{\lam}$ has at least two solutions for $\lam \in ( \lam _2^*,
\lam ^*)$.
    \end{enumerate}
\end{thm}
A uniqueness result  in the spirit of (1) also holds for dimension
$1$ (resp., dimension $2$) with $N/2$ replaced by $1$ (resp.,
$1+\epsilon$). However, in spite of above results, issues of
uniqueness, multiplicity and other qualitative properties of the
solutions for $(S)_{ \lam}$ are still far from being well
understood. For example, we conjecture that no solution exists for
$(S)_{\lambda^*}$ with $N\geq 8$ --at least when $f\equiv 1$. In \S
5 we shall present some numerical evidences for various conjectures
relating to the case of power-law permittivity profile
$f(x)=|x|^\alpha$ defined in a unit ball. It looks like there are
two critical exponents $\alp ^*=-\frac{1}{2}+\frac{1}{2}\sqrt{27/2}$
and $ \alp ^{**}(N)=\frac{4-6N+3\sqrt{6}(N-2)}{4}$ (which is
relevant for $N\ge 8$) such that the following four regimes are
possible:
\begin{enumerate}
\item  There exist exactly two solutions for $0<\lambda<\lambda^*$, and
one solution for $\lambda=\lambda^*$. This regime occurs when
   $N=1$ and $\alpha \leq 1$.

\item There exists exactly one solution for $0<\lambda < \lambda^*_1$,
exactly 2 solutions for $\lambda^*_1<\lambda <\lambda^*$ and exactly
one at $\lambda=\lambda^*$. This regime occurs when
   $N=1$ and $1\leq  \alpha \leq \alpha^*$.
   \item There exists exactly one solution for $0<\lambda <\lambda^*_1$,
exactly two solutions for $\lambda_2^* <\lambda <\lambda^*$, while
multiple solutions can be obtained for $\lambda_1^* <\lambda
<\lambda_2^*$. Moreover, the multiplicity becomes arbitrarily large
as $\lambda$ approaches another critical value $\lambda_*\in
(\lambda_1^*, \lambda_2^*)$, at which there is a touchdown
(quenching) solution $u$ characterized with $\parallel u
\parallel_{\infty}=1$. This regime occurs when
\begin{itemize}
\item $2\leq N\leq 7$  and $\alpha \geq 0$;
\item $N\geq 8$ and $\ \alpha^{**}<\alpha$.
\end{itemize}
\item There is exactly one solution if  $0<\lambda<\lambda^*$ and none for
$\lambda \geq \lambda^*$. This regime occurs when
   $N\geq 8$,  and $0\leq \alpha \leq \alpha^{**}$.

\end{enumerate}
We finally mention that the above results  can be extended to more
general elliptic problems of the form
$$
\arraycolsep=1.5pt\begin{array}{lll}\arraycolsep=1.5pt
     \quad\quad   -\Delta v &=& \displaystyle\frac{\lambda
f(x)}{(1-v)^{\beta}} \,,\  \ x \in \Omega \,;\\
    \quad\quad\  v (x)  &=& 0  ,\quad\quad\quad\quad x \in \partial\Omega
\,
\end{array}\eqno{(S)_{\lam, \beta}}$$
with $\beta >0$. Here the critical dimension depends on the
parameter $\beta $, and this is the subject of a work in progress.

\section{Pull-In Voltage}

In this section, we study the steady-state deflection $u$ which
satisfies $(S)_\lambda$, and we establish the existence and some
estimates on the pull-in voltage $\lam ^*$ for $(S)_{ \lam }$
defined as:
\begin{equation}
\lam ^*=\sup\{\lam >0\ |\ (S)_{ \lam } \ possesses \ at \ least\
one\ solution\} \,.\label{1:lam}
\end{equation}
In other words, $\lam ^*$ is called pull-in voltage if there exist
uncollapsed states for $0<\lam <\lam ^*$ while there are none of
them for $\lam >\lam ^*$.

\begin{thm} There exists a finite pull-in voltage $\lam ^*>0$ such that
\begin{enumerate}
\item If $\lam <\lam ^*$, there exists at least one solution for
$(S)_{ \lam}$;
\item If  $\lam >\lam ^*$, there is no solution for $(S)_{\lam}$.
\end{enumerate}
Moreover, with $\nu _{_\Omega} $ defined by (1.6), we have the lower
bound
\begin{equation}
\nu _{_\Omega}  \Big(\sup_{x\in \Omega} f(x)\Big)^{-1}\le\lam
^*\,.\label{0:6}
\end{equation}

\end{thm}

\noindent{\bf Proof:} We need to show that $(S)_{ \lam }$ has at
least one solution when $ \lam < \nu_{_\Omega} (\sup_{\Omega}
f(x))^{-1}. $ Indeed, it is clear that $u\equiv 0$ is a sub-solution
of $(S)_{ \lam }$ for all $\lam >0$. To construct a super-solution
of $(S)_{ \lam }$, we consider a bounded domain $\Gamma \supset \bar
\Omega$ with smooth boundary, and let $(\mu_{_\Gamma},
\psi_{_\Gamma})$ be its first eigenpair normalized in such a way
that
\[
\hbox{$\sup\limits_{x\in \Gamma} \psi_\Gamma(x)=1$ and
$\inf\limits_{x\in \Omega} \psi_\Gamma(x):=s_1 >0.$}
\]
   We
construct a super-solution in the form $\psi=A\psi _{_\Gamma}$ where
$A$ is a scalar to be chosen later. First, we must have $A\psi
_{_\Gamma}\ge 0$ on $\partial \Omega $ and $0<1-A\psi_{{_\Gamma}}<
1$ in $\Omega $, which requires that
\begin{equation}0<a<1. \label{1:7}
\end{equation}
We also require
\begin{equation}
-\Delta  \psi - \frac{\lam  f(x)} {(1-A \psi)^2}\geq 0 \ \quad {\rm
in} \quad \Omega \,,
\end{equation}
which can be satisfied as long as:
\begin{equation} \mu _{{_\Gamma}}A\ \psi _{_\Gamma}\ge\frac{\lam
\sup_{\Omega} f(x)} {(1-A\ \psi _{_\Gamma})^2}\ \quad {\rm in} \quad
\Omega \,,
\end{equation}
or
\begin{equation}
\lambda\sup\limits_\Omega f(x) < \beta (A,
\Gamma):=\mu_{_\Gamma}\inf \{g(sA); \, s\in [s_1(\Gamma),1]\}\,,
\end{equation}
where $g(s)=s(1-s)^2$. In other words, $\lambda^*\sup\limits_\Omega
f(x)\geq \sup\{ \beta (A, \Gamma); \, 0<a<1, \Gamma \supset \bar
\Omega\}$, and therefore  it remains to show that
\begin{equation}
\nu_\Omega=\sup\{ \beta (A, \Gamma); \, 0<a<1, \Gamma \supset \bar
\Omega\}.
\end{equation}
For that, we note first that
\[
\inf\limits_{s\in [s_1,1]}g(As)=\min \big\{g(As_1), g(A)\big\}.
\]
   We also have that
$g(As_1)\leq g(A)$ if and only if $A^2(s_1^3-1)-2A(s_1^2-1)+(s_1-1)
\leq 0$ which happens if and only if $A^2(s_1^2+s_1+1)-2A(s_1+1)+1
\geq 0$ or if and only if either $A\leq A_-$ or $A\geq A_+$ where
$$A_+=\frac{s_1+1+\sqrt{s_1}}{s_1^2+1+s_1}=\frac{1}{s_1+1-\sqrt{s_1}}\,,
\quad
A_-=\frac{s_1+1-\sqrt{s_1}}{s_1^2+1+s_1}=\frac{1}{s_1+1+\sqrt{s_1}}
$$
Since $A_-<1<a_+$, we get that
\begin{equation}
G(A)=\inf\limits_{s\in
[s_1,1]}g(As)=\left\{\begin{array}{lll}g(As_1)
\quad &{\rm if}&\quad 0\leq A \leq A_-\,,\\[3mm]
g(A) \quad &{\rm if}&\quad A_-\leq A\leq 1\,.
\end{array}\right.
\end{equation}
We now have  that $\frac{dG}{dA}=g'(As_1)s_1\geq 0$ for all $0\leq A
\leq A_-$.  And since $A_-\geq \frac{1}{3}$, we have
$\frac{dG}{dA}=g'(A)\leq 0$ for all $A_-\leq A \leq 1$. It follows
that $$\arraycolsep=1.5pt \begin {array}{lll} \arraycolsep=1.5pt
\sup\limits_{0<a<1}\inf\limits_{s\in [s_1,1]}
g(As)&=&\sup\limits_{0<a<1}G(A)=G(A_-)=g(A_-)\\[3mm]
&=&\displaystyle \frac{1}{s_1+1+\sqrt{s_1}}\big(1-\displaystyle
\frac{1}{s_1+1+\sqrt{s_1}}\big)^2\\[3mm]
&=&\displaystyle
\frac{s_1(s_1+1+2\sqrt{s_1})}{(s_1+1+\sqrt{s_1})^3}\\[3mm]
&=&H(\inf_\Omega \psi_\Gamma))
\end{array}$$
which proves our lower estimate.

Now that we know that $\lambda^*>0$, pick $\lam \in (0, \lam ^*)$
and use the definition of $\lambda^*$ to find a $\bar {\lam}\in
(\lam ,\lam ^*)$ such that $(S)_{ \bar {\lam}}$ has a solution
$u_{\bar {\lam}}$, $i.e$,
$$
     -\Delta u _{\bar {\lam}}= \frac{\bar {\lam } f(x)}{(1-u_{\bar
{\lam}})^2} \,, \quad x \in \Omega \,;
      \qquad u_{\bar {\lam}}=0 \,, \quad x \in \partial\Omega \,.
$$
and in particular $ -\Delta u _{\bar {\lam}}\ge  \frac{\lam
f(x)}{(1-u_{\bar {\lam}})^2}$ for $ x \in \Omega $ which then
implies that $u_{\bar {\lam}}$ is a super-solution of $(S)_{ \lam
}$. Since $u\equiv 0$ is a sub-solution of $(S)_{ \lam }$, then we
can again conclude that there is a solution $u_{\lam }$ of $(S)_{
\lam }$ for every $\lam \in (0, \lam ^*)$.

It is also easy to show that $\lambda^*$ is finite, since if $(S)_{
\lam }$ has at least one solution $0<u<1$, then by integrating
against the first (positive) eigenfunction $\phi_{_\Omega}$, we get
\begin{equation}
+\infty >\mu_{_\Omega}\geq \mu_{_\Omega} \displaystyle \int _{\Omega
}u\phi _{_\Omega}=-\displaystyle \int _{\Omega }u\Delta
\phi_{_\Omega} =-\displaystyle \int _{\Omega }\phi_{_\Omega}\Delta
u=\lam \displaystyle \int _{\Omega }\frac{\phi_{_\Omega}
f}{(1-u)^2}dx \geq \lam \displaystyle \int _{\Omega }\phi_{_\Omega}
f dx  \label{1:5i}
\end{equation}
and therefore $\lam ^*<+\infty $. The definition of $\lam ^*$
implies that there is no solution of $(S)_{\lam }$ for any $\lam
>\lam ^*$.
\endproof

\subsection{Lower bounds for $\lam^*$}

   It is  desirable to seek more computationally accessible lower bounds on
$\lam ^*$. For that we consider for every subset $\Gamma \subset
\mathbb{R}^N$ and any function $f$ on $\Gamma$ such that $0\leq f
\leq 1$, the corresponding  pull-in voltage $\lam ^*(\Gamma, f)$,
that is the value $\lambda^*$ defined above for  the problem \bsub
\label{2:31} \arraycolsep=1.5pt\arraycolsep=1.5pt
\begin{eqnarray}
     -\Delta u &=&\frac{\lambda f(x)}{(1-u)^{2}} \ \  \quad   \quad x \in
\Gamma \,, \label{i:s1}\\
\hfill 0<&u&<1 \quad \quad \quad \quad \ \ x \in
\Gamma \,, \label{i:s2}\\
u&=&0  \quad \quad \quad \qquad     x \in
\partial\Gamma .
\label{i:s3}
\end{eqnarray} \esub

  We need the following result  which can be found in \cite{B} (Theorem
4.1).

\begin{lem} For any bounded domain $\Gamma$ in $\R^N$ and any function $f$
on $\Gamma$ such that $0\leq f \leq 1$, we have
\[
\lam ^*(\Gamma, f)\geq \lam ^*(B_R, f^*)
\]
   where $B_R=B_R(0)$ is the Euclidean ball in $\mathbb{R}^N$ with
radius $R>0$ and with volume $ |B_R|=|\Gamma |$, and where $f^*$ is
the Schwarz symmetrization of $f$.
\end{lem}

We now establish the following refined lower bounds for $\lam ^* $
of $(S)_{\lam}$.

\begin{lem}
We have the following lower bound for $\lambda^*$.
\begin{equation}\lam ^*\ge \frac{8N}{27\sup_{\Omega}f}\Big(\frac{\omega
_{_N} }{|\Omega|}\Big)^{\frac{2}{N}}. \label{2:32}
\end{equation}
Moreover, if $f(x)\equiv |x|^{\alp}$ on $\Omega$ with $\alp \ge 0$,
then we have
\begin{equation}\lam ^*\ge \frac{4(2+\alp)(N+\alp) }{27}\Big(\frac{\omega
_{_N} }{|\Omega|}\Big)^{\frac{2+\alp}{N}}. \label{2:33}
\end{equation}
\end{lem}

\noindent  {\bf Proof:} Setting $R=\Big(\frac{|\Omega| }{\omega
_{_N}}\Big)^{\frac{1}{N}}$, it suffices  --in view of Lemma 2.2--
and since $\sup\limits_{B_R}f^*=\sup\limits_{\Omega}f$, to show that
\begin{equation}\lam ^*\ge \frac{8N }{27R^2\sup_{\Omega}f}\,,
\label{2:34}
\end{equation}
for the case where $\Omega =B_R$. In fact, the function
$w(x)=\frac{1}{3}(1-\frac{|x|^2}{R^2})$  satisfies on $B_R$
$$ \arraycolsep=1.5pt
\begin{array}{lll}
-\Delta w= \displaystyle\frac{2N}{3R^2}&=&
\displaystyle\frac{2N(1-\frac{1}{3})^{2}}{3R^2}
\displaystyle\frac{1}{(1-\frac{1}{3})^{2}}\\[3mm]
&\ge  &\displaystyle\frac{8N}{27R^2\sup_{\Omega}f}\displaystyle
\frac{f(x)}{[1-\frac{1}{3}(1-\frac{|x|^2}{R^2})]^{2}}\\[3mm]
&=&\displaystyle\frac{8N}{27R^2\sup_{\Omega}f}\displaystyle\frac{f(x)}{(1-w)^{2}}\,.
\end{array}
$$
So for $ \lam \le \frac{8N}{27R^2\sup_{\Omega}f}\,, $ $w$ is a
super-solution of $(S)_{\lam}$ in $B_R$. Since on the other hand
$w_0\equiv 0$ is a sub-solution of $(S)_{\lam}$ and $w_0\le w$ in
$B_R$, then there exists  a solution of $(S)_{\lam}$ in $B_R$ which
proves $(\ref{2:34})$ and hence $(\ref{2:32})$.

In order to prove $(\ref{2:33})$, it suffices to note that
    $w(x)=\frac{1}{3}\big(1-\frac{|x|^{2+\alp}}{R^{2+\alp}}\big)$  is a
super- solution for  $(S)_{\lam}$
   on $B_R$ provided  $ \lam
\le \frac{4(2+\alp)(N+\alp)}{27R^{2+\alp}}$. This completes the
proof of Lemma 2.3.
\endproof

\subsection{Upper bounds for $\lam^*$}

We note that $(\ref{1:5i})$ already yields a finite upper bound for
$\lam ^*$. However, Pohozaev-type arguments such as the one used in
\cite{GPW} can be used to establish better and more computable upper
bounds. In this subsection, we establish these estimates and hence
complete the proof of Theorem 1.1.

We shall consider problem $(S)_{\lam }$ in the case where $\Omega
\subset\mathbb{R}^N$  is a strictly star-shaped domain containing
$0$, meaning that $\Omega $ satisfies the additional property that
there exists a positive constant $a$ such that
\begin{equation}
x\cdot \nu \ge a\int _{\partial \Omega }ds \quad  {\rm for\ all}
\quad x\in\partial \Omega \,, \label{2:1i}
\end{equation} where $\nu $
is the unit outer normal to $\Omega $ at $x\in\partial \Omega $.

\begin{prop} Suppose $f\equiv 1$ and that the strictly star-shaped domain
$\Omega \subset\mathbb{R}^N$ satisfies $(\ref{2:1i})$. Then the
pull-in voltage $\lambda^*(\Omega)$ satisfies:
\begin{equation}
\lam ^*(\Omega)\le \bar{\lambda} _3=\frac{(N+2)^2}{8aN|\Omega |}\,.
\label{2:2i}
\end{equation}
In particular, if  $\Omega =B_1(0)\subset\mathbb{R}^N$ then we have
the bound
\[
\lam ^*(B_1(0))\le \frac{(N+2)^2}{8 }\,.
\]
\end{prop}

\noindent  {\bf Proof:} Recall the well-known Pohozaev's identity:
If  $u$ is a solution of
\[
\begin{array}{lll}
\Delta u+\lam g(u)=0 \quad \quad{\rm for}
      \quad x \in \Omega  \,,\\
     \quad \quad \quad \quad \,  \ u=0    \quad \quad  {\rm for}
\quad x \in \partial\Omega \,,
\end{array}
\]
then
\begin{equation}
N\lam \int _{\Omega }G(u)dx-\frac{N-2}{2}\lam \int _{\Omega
}ug(u)dx=\frac{1}{2}\int _{\partial\Omega }(x\cdot\nu
)\big(\frac{\partial u}{\partial \nu}\big)^2ds\,,\label{2:3i}
\end{equation}
where $G(u)=\int ^{u}_{0}g(s)ds$. Applying it with $
g(u)=\frac{1}{(1-u)^2}$ and $G(u)=\frac{u}{1-u}$, it yields
\begin{equation}\arraycolsep=1.5pt \begin {array}{lll}
    \displaystyle\frac{\lambda}{2}\displaystyle\int _{\Omega }\displaystyle
\frac{u(N+2-2Nu)}{(1-u)^2}dx&=&
\displaystyle\frac{1}{2}\displaystyle\int _{\partial\Omega }(x\cdot
\nu )\big(\displaystyle\frac{\partial u}{\partial
\nu}\big)^2ds\\[4mm]
&\ge &\displaystyle\frac{a}{2}\Big(\displaystyle\int
_{\partial\Omega }\displaystyle\frac{\partial u}{\partial
\nu}ds\Big)^2\\[4mm]
&=&\displaystyle\frac{a}{2}\Big(-\displaystyle\int _{\Omega
}\Delta udx\Big)^2\\[4mm]
&=&\displaystyle\frac{a\lam ^2}{2}\Big(\displaystyle\int _{\Omega }
\displaystyle\frac{dx}{(1-u)^2}\Big)^2\,,
\end{array} \label{2:4i}\end{equation}
where we have used the Divergence Theorem and Holder's inequality
\[
\int _{\partial\Omega }\frac{\partial u}{\partial \nu}ds\le
\Big(\int _{\partial\Omega }\big(-\frac{\partial u}{\partial
\nu}\big)^2ds\Big)^{1/2}\Big(\int _{\partial\Omega }ds\Big)^{1/2}\,.
\]
Since
\[\arraycolsep=1.5pt \begin {array}{lll}
\displaystyle \int _{\Omega }\displaystyle
\frac{u(N+2-2Nu)}{(1-u)^2}dx&=&\displaystyle \int _{\Omega
}\Big[-2N\big(u-\displaystyle\frac{N+2}{4N}\big)^2+
\displaystyle\frac{(N+2)^2}{8N}\Big]\displaystyle
\frac{1}{(1-u)^2}dx\\[4mm]
&\le &\displaystyle\frac{(N+2)^2}{8N}\displaystyle\int _{\Omega
}\frac{dx}{(1-u)^2}\,,
\end{array}\]
we deduce from $(\ref{2:4i})$ that
\[
\frac{(N+2)^2}{8N}\ge a\lam \int _{\Omega }\frac{dx}{(1-u)^2}\ge
a\lam |\Omega |\,,
\]
which implies the upper bound $(\ref{2:2i})$ for $\lam ^*$.

Finally, for the special case where $\Omega
=B_1(0)\subset\mathbb{R}^N$, we have $a=\frac{1}{N\omega _{_N}}$
with $\omega _{_N}=|B_1(0)|$, and hence the bound $\lam ^*\le
\bar{\lambda}_3=\frac{(N+2)^2}{8 }$.
\endproof

For a general domain $\Omega$, the following upper bounds on
$\lam^*(\Omega)$ established in \cite{P1} and  \cite{GPW}
respectively, complete the proof of Theorem 1.1.

\begin{prop} (1) If  $f$ satisfies $ 0<C_0\le f(x) \le 1$ on $\Omega
$, then \begin{equation}
    \lam ^*\le \bar{\lam}_1 \equiv \frac{4\mu _{_\Omega}
}{27C_0} \,.\label{p:lam1}
\end{equation}

(2) If $f$ satisfies $0 \le f(x) \le 1$ on $ \Omega$, and  if $f>0$
on a set of positive measure, then
     \begin{equation}
     \lam_{*} \le \bar{\lam}_2 \equiv \frac{\mu _{_\Omega}}{3} \left(
\int_{\Omega}
      f \phi_{_\Omega} \, dx \right)^{-1} \,.\label{p:lam2}
\end{equation}
Here $\mu_{_\Omega} $ and $\phi_{_\Omega}$ are the first eigenpair
of $-\Delta$ on $H^1_0(\Omega)$ with $\int _{\Omega
}\phi_{_\Omega}dx=1$.
\end{prop}

\subsection{Numerical estimates for $\lam^*$}

\begin{table}
\begin{center}  Exponential Profiles: \\[2mm]
\begin{tabular}{  c | c | c | c|c | c  }
\hline $\Omega$ & $\alpha$ &$\underline\lam $& $\lam^{*}$ &
$\bar{\lam}_1$ & $\bar{\lam}_2$
\\
\hline
\mbox{(Slab)}            & 0 &    1.185  &1.401  & 1.462  & 3.290  \\
    \mbox{(Slab)}           &  1.0 &  1.185  &1.733  & 1.878  & 4.023  \\
    \mbox{(Slab)}           &  3.0 &  1.185  &2.637  & 3.095  & 5.965  \\
    \mbox{(Slab)}           &  6.0 &  1.185  &4.848  & 6.553  & 10.50  \\
    \mbox{(Unit Disk)}      & 0 &    0.593  &0.789  & 0.857 & 1.928  \\
    \mbox{(Unit Disk)}      &  0.5 &  0.593 &1.153  & 1.413  & 2.706  \\
    \mbox{(Unit Disk)}      &  1.0 &  0.593 &1.661  & 2.329  & 3.746  \\
    \mbox{(Unit Disk)}      &  3.0 &  0.593 &6.091  & 17.21  & 11.86  \\
\hline \end{tabular}
\end{center}
\caption{{\em Numerical values for pull-in voltage $\lam ^{*}$ with
the bounds given in Theorem 1.1. Here the exponential permittivity
profile is chosen as (\ref{p:perm}).}}
\end{table}

     In the computations below we shall consider two choices for the domain
$\Omega$,
\begin{equation}
     \Omega: [-{1/2},{1/2}] \quad \mbox{(Slab)} \,; \quad
     \Omega: x^2 + y^2 \leq 1 \quad \mbox{(Unit Disk)} \,. \quad
      \label{p:dom}
\end{equation}
To compute the bounds $\bar{\lam}_1$ and $\bar{\lam}_2$, we must
calculate the first eigenpair $\mu_{_\Omega}$ and $\phi_{_\Omega}$
of $-\Delta$ on $\Omega$, normalized by $\int_{\Omega}
\phi_{_\Omega}\, dx =1$, for each of these domains.  A simple
calculation yields that \bsub \label{p:eig_ex}
\begin{align}
      \mu_{_\Omega}&= \pi^2 \,, \quad \quad\qquad \phi_{_\Omega} =
\frac{\pi}{2} \sin\left[\pi \left(
      x + \frac{1}{2} \right) \right] \,, \quad \mbox{(Slab)} \,;
         \label{p:eig_exa}  \\
      \mu_{_\Omega} & = z_0^{2}\approx 5.783\,, \qquad \phi_{_\Omega} =
\frac{z_0}{J_{1}(z_0)} J_{0}(z_0 |x|) \,, \quad \mbox{(Unit Disk)}
\,. \label{p:eig_exb}
\end{align}
\esub Here $J_0$ and $J_1$ are Bessel functions of the first kind,
and $z_0\approx 2.4048$ is the first zero of $J_0(z)$. The bounds
$\bar{\lam}_1$ and $\bar{\lam}_2$ can be evaluated by substituting
(\ref{p:eig_ex}) into (\ref{p:lam1}) and (\ref{p:lam2}). Notice that
$\bar{\lam}_2$ is, in general, determined only up to a numerical
quadrature.

Using Newton's method and COLSYS \cite{A}, one can also solve the
boundary value problem $(S)_{\lam }$ and numerically calculate
$\lam^{*}$ as the saddle-node point for the following two choices of
the permittivity profile: \bsub \label{p:perm}
\begin{align}
    \mbox{(Slab)}: \quad & f(x) = |2x|^{\alp}\,, \quad \mbox{(power-law)}
\,;\quad
    f(x)=e^{\alp (x^2-{1/4})}  \quad \mbox{(exponential)}\,,
\label{p:perma}\\
\mbox{(Unit Disk)}: \quad & f(x) = |x|^{\alp} \,, \quad \, \
\mbox{(power-law)} \,;
    \quad f(x)=e^{\alp (|x|^2-1)}  \,, \quad \mbox{(exponential)}\,,
     \label{p:permb}
\end{align}
\esub where $\alp\ge 0$.  Table 1 contains numerical values for
$\lam ^*$ in the case of  exponential profiles, while  Table 2 deals
with  power-law profiles.
\begin{table}
\begin{center}   Power-Law Profiles:\\[2mm]
\begin{tabular}{c|c|c|c|c|c}
\hline
    $\Omega$ & $\alpha$ & $\lam _c(\alp)$&$\lam ^{*}$ & $
\bar{\lam}_1$ &  $\bar{\lam}_2$
\\
\hline
    \mbox{(Slab)}           &  0& 1.185  & 1.401 &  1.462  & 3.290   \\
    \mbox{(Slab)}           &   1.0&  3.556 &4.388  &  $\infty$  &  9.044
\\
    \mbox{(Slab)}           &  3.0 &  11.851 &15.189  &  $\infty$   &
28.247
\\
    \mbox{(Slab)}           &  6.0&  33.185 &43.087 & $\infty$    & 76.608
\\
    \mbox{(Unit Disk)}      &  0  &  0.593 &0.789 &  0.857  &  1.928  \\
    \mbox{(Unit Disk)}      &  1.0 & 1.333  &1.775 & $\infty$   &  3.019
\\
    \mbox{(Unit Disk)}      &  5.0&  7.259 &9.676  &  $\infty$   &   15.82
\\
    \mbox{(Unit Disk)}      &  20&  71.70 &95.66 & $\infty$   &   161.54 \\
\hline
\end{tabular}
\end{center}
\caption{{\em Numerical values for pull-in voltage $\lam ^{*}$ with
the bounds given in Theorem 1.1. Here the power-law permittivity
profile is chosen as (\ref{p:perm}).}}
\end{table}
    What is remarkable is that $\bar{\lam}_1$ and
$\bar{\lam}_2$ are not comparable even when $f$ is bounded away from
$0$ and that neither one of them provides the optimal value for
$\lambda^*$. This leads us to conjecture that there should be a
better estimate for $\lambda^*$, one  involving the distribution of
$f$ in $\Omega$, as opposed to the infimum or its average against
the first eigenfunction $\phi_{_\Omega}$.

\section{Minimal Positive Solutions}
   In this section, we are
concerned with {\it minimal positive solutions} for $(S)_{\lam }$.
We establish their existence, uniqueness and other related
properties.  We consider the case $\lam \in (0,\lam ^*)$ in \S 3.1,
and   $\lam =\lam ^*$ in \S 3.3, but first we give a recursive
scheme for the construction of such solutions.

\begin{thm}
  For any $\lam <\lam^*$ there exists a unique minimal positive
solution $u_{\lam}$ for $(S)_{\lam }$. It is obtained as the limit
of the sequence $\{u_n(\lam ; x)\}$ constructed recursively as
follows: $u_0\equiv 0$ in $\Omega $ and for each $n\geq 1$,
\begin{equation}\arraycolsep=1.5pt \begin {array}{lll}
-\Delta u_n &= \displaystyle\frac{\lambda f(x)}{(1-u_{n-1})^{2}} \,,
\quad &x \in \Omega \,;\\[3mm]
\quad 0\le &u_n<1 \,, \quad \quad\quad &x \in \Omega \,;
\\[3mm]
\quad \  u_n &= 0 \,,\quad \quad \quad \quad &x \in
\partial\Omega.
\end{array} \label{3:1}\end{equation}
\end{thm}

\noindent  {\bf Proof:} Let $u$ be any positive solution for
$(S)_{\lam }$, and consider the sequence $\{u_n(\lam ; x)\}$ defined
in $(\ref{3:1})$. Clearly $u( x)>u_0\equiv 0$ in $\Omega $, and
whenever $u(x)\ge u_{n-1}$ in $\Omega $, then
$$\begin {array}{lll}
      -\Delta (u-u_n) = \lambda f(x)\big [\displaystyle\frac{1}{(1-u)^{2}}-
\displaystyle
      \frac{1}{(1-u_{n-1})^{2}}\big ]\ge 0
      \,, \quad x \in \Omega \,\\[3mm]
     \quad \ \ \ \   u-u_n=0 \,,\ \ \ \ \quad x \in \partial\Omega \,.
\end{array} $$
The maximum principle and an immediate induction yield that
$1>u(x)\ge u_{n}$ in $\Omega $ for all $n\ge 0$. In a similar way,
the maximum principle implies that the sequence $\{u_n(\lam ; x)\}$
is monotone increasing. Therefore, $\{u_n(\lam ; x)\}$ converges
uniformly to a positive solution $u_{\lam}(x)$, satisfying $u( x)\ge
u_{\lam}(x)$ in $\Omega $, which is a minimal positive solution of
$(S)_{\lam }$. It is also clear that $u_{\lam}(x)$ is unique.
\endproof

\begin{rem} Let $g(x,\xi, \Omega )$ be the Green's function of Laplace
operator, with $g(x,\xi, \Omega )=0$ on $\partial\Omega$. Then the
iteration in  $(\ref{3:1})$  can be replaced by: $u_0\equiv 0$ in
$\Omega $ and for each $n\geq 1$,
\begin{equation}\arraycolsep=1.5pt \begin {array}{lll} \arraycolsep=1.5pt
\quad \ u_n(\lam;x)&=\lam \displaystyle\int
_{\Omega}\displaystyle\frac{f(\xi)g(x,\xi, \Omega)}{(1-u_{n-1}(\lam
;\xi))^{2}}d\xi \,,\quad &x \in \Omega \,;\\
\quad \  u_n (\lam;x)&=0 \,,\quad \quad \quad \quad &x \in
\partial\Omega\,.
\end{array} \label{3:2+}\end{equation}
The same reasoning as above yields that
$\lim_{n\to\infty}u_n(\lam;x)=u_{\lam }(x)$ for all $ x\in\Omega$.
\end{rem}
The above construction of solutions yields the following
monotonicity result for the pull-in voltage.
\begin{prop}
If  $\Omega _1\subset\Omega _2$, then   $\lam ^*(\Omega _1)\ge
\lam^*(\Omega _2)$ and the corresponding minimal solutions satisfy
$u_{{_{\Omega _1}}}(\lambda, x)\le  u_{{_{\Omega _2}}}(\lambda, x)$
on $\Omega _1$  for every $0<\lambda <\lam^*(\Omega _2)$.
\end{prop}

\noindent  {\bf Proof:} Again the method of sub/super-solutions
immediatly yields that $\lam ^*(\Omega _1)\ge \lam ^*(\Omega _2)$.
Now consider for $i=1, 2$, the sequences $\{u_n(\lambda, x,\Omega
_i)\}$ on $\Omega _i$ defined by $(\ref{3:2+})$ where
$g(x,\xi,\Omega _i)$ are the corresponding Green's functions on
$\Omega _i$. Since $\Omega _1\subset\Omega _2$, we have that
$g(x,\xi,\Omega _1)\le g(x,\xi,\Omega _2)$ on $\Omega _1$. Hence, it
follows that
\[
u_1(\lambda, x,\Omega _2)=\lam \int_{\Omega _2}f(\xi)g(x,\xi,\Omega
_2)d\xi\ge \lam \int_{\Omega _1}f(\xi)g(x,\xi,\Omega
_1)d\xi=u_1(\lambda, x,\Omega _1)
\]
on $\Omega _1$. By induction we conclude that $u_n(\lambda, x,\Omega
_2)\ge u_n(\lambda, x,\Omega _1)$ on $\Omega _1$ for all $n$. On the
other hand, since
   $u_n(\lambda, x,\Omega _2)\le u_{n+1}(\lambda, x,\Omega _2)$ on $\Omega
_2$ for $n$, we get that $u_n(\lambda, x,\Omega _1)\le u_{_{\Omega
_2}}(\lambda, x)$ on $\Omega _1$, and we are done.
\endproof

\subsection{Spectral properties of minimal solutions}

For a further study of  minimal (positive) solutions, we now
consider for each positive solution $u$ of $(S)_\lambda$, the
operator
\begin{equation}
L_{u, \lam }=-\Delta -\frac{2\lam f}{(1-u)^{3}}
\end{equation}
associated to the linearized problem around $u$. We see that minimal
solutions in the above sense correspond to variational solutions
that are local minimizers. We denote by $\mu (\lam, u )$ the
smallest eigenvalue of $L_{u, \lam }$, that is the one corresponding
to the following Dirichlet eigenvalue problem \bsub
    \arraycolsep=1.5pt\arraycolsep=1.5pt
\begin{eqnarray}
-\Delta \phi -\frac{2\lambda f(x)}{(1-u)^{3}}\phi &=&\mu (\lam, u
)\phi\,, \quad    x \in \Omega \,;\\
\phi &=&0  \quad \quad  \quad \ \   x \in
\partial\Omega \,.
\end{eqnarray}
\esub In other words,
$$
\mu (\lam, u )=\inf _{\phi \in H^1_0(\Omega )}\frac{\int _{\Omega
}\big \{|\nabla \phi |^2-2\lam f(1- u)^{-3}\phi ^2\big\}dx}{\int
_{\Omega }\phi ^2dx}\,.
   $$
\begin{prop} The following hold:
\begin{enumerate}
\item
$\lam ^{*}=\sup\{\lam ;  \hbox{ $L_{u _{_{\lam}}, \lam }$ has
positive first eigenvalue for minimal solution $u _{\lam}$ of
$(S)_\lambda$} \}.  $
\item If $0<\lam <\lam ^*$  then the smallest eigenvalue
$\mu_\lambda:=\mu (\lambda, u_\lambda)$ of $L_{u_\lambda, \lam }$
--corresponding to the minimal solution $u_\lambda$-- is positive
and $\lambda \to \mu_\lambda$ is decreasing on $(0, \lambda^*)$.
\end{enumerate}
    \end{prop}

For Proposition 3.3, we need the following crucial lemma.

\begin{lem}
Suppose $u$ is a positive solution of $(S)_{\lam}$, and let $\mu
(\lambda, u)$ be the corresponding first eigenvalue. Consider any
-classical- supersolution $v$ of $(S)_\lambda$, that is
    \bsub \label{3:2}
\arraycolsep=1.5pt\arraycolsep=1.5pt
    \begin{eqnarray}
      -\Delta v &\ge &\frac{\lambda f(x)}{(1-v)^{2}} \ \  \quad   \quad x
\in
\Omega \,, \label{i:s1}\\
      0\le& v(x)&<1 \quad \quad \quad x\in \Omega \\
    v&=&0  \quad \quad \quad \qquad    \quad x \in
\partial\Omega.
\label{i:s3}\end{eqnarray} \esub \ If $\mu(\lambda, u)>0$ then $v\ge
u$ on $\Omega $, and if $\mu (\lambda, u)=0$ then $v= u$ on $\Omega
$.

\end{lem}
\noindent  {\bf Proof:} For a given $\lam $ and $x\in \Omega$, use
the fact that $f(x)\ge 0$  and that $t\to \frac{\lambda
f(x)}{(1-t)^{2}}$ is convex on $(0,1)$, to obtain
\begin{equation}
-\Delta (u+\tau (v-u))-\frac{\lam f(x)}{[1-(u+\tau (v-u))]^{2 }}\ge
0 \quad x\in \Omega\,, \label{3:4}
\end{equation}
for $\tau \in [0,1]$. Note that $(\ref{3:4})$ is an identity at
$\tau =0$, which means that the first derivative of the left side
for $(\ref{3:4})$ with respect to $\tau$ is nonnegative at $\tau
=0$, $i.e$,
    \bsub \label{3:5}
\arraycolsep=1.5pt\arraycolsep=1.5pt   \begin{eqnarray}
      -\Delta (v-u) -\frac{2\lambda f(x)}{(1-u)^{3}}(v-u)&\ge &0 \
      \quad   \quad x \in \Omega \,, \label{i:s1}\\
    v-u&=&0  \quad  \,    \quad x \in
\partial\Omega \,.
\label{i:s3}\end{eqnarray} \esub Thus, the maximal principle implies
that if $\mu (\lam ,u)>0$ we have $v\ge u$ on $\Omega $, while if
$\mu (u)=0$ we have\begin{equation}
    -\Delta (v-u) -\frac{2\lambda f(x)}{(1-u)^{3}}(v-u)=0 \   \quad   \quad
x \in \Omega \,. \label{3:6}
\end{equation}
In the latter case the second derivative of the left side for
$(\ref{3:4})$ with respect to $\tau$ is nonnegative a $\tau =0$
again, $i.e$,
\begin{equation}
-\frac{6\lam f(x)}{(1-u)^{4}}(v-u)^2\ge 0 \quad x\in \Omega\,,
\label{3:7}
\end{equation}
  From $(\ref{3:7})$ we deduce that $v\equiv u$ in $\Omega \setminus
\Omega _0 $, where
\begin{equation}
\Omega _0 =\{x\in \Omega :\, f(x)=0\ \mbox{for}\  x\in \Omega \}\,.
\label{3:8}
\end{equation}
On the other hand,   $(\ref{3:6})$ reduces to
$$\begin{array}{ll}
    -\Delta (v-u) =0\quad \ x\in \Omega _0 \,,\\
\quad \quad \ \,  v-u =0 \ \quad x\in \partial\Omega _0\,,
\end{array}
$$
which implies $v\equiv u$ on $\Omega _0$. Hence if $\mu (\lam ,u)=0$
then $v\equiv u$ on $\Omega $, which completes the proof of Lemma
3.4.
\endproof

\noindent  {\bf Proof of Proposition 3.3:} (1) Let
\[
\lam ^{**}=\sup\{\lam ;  \hbox{ $L_{u _{_{\lam}}, \lam }$ has
positive first eigenvalue for minimal solution $u _{\lam}$ of
$(S)_\lambda$} \}.
\]
It is clear that $\lam ^{**}\le \lam ^*$, so it suffices to prove
that there is no minimal solution for $(S)_{\mu}$ with $\mu >\lam
^{**}$. In fact, suppose $w$ is a minimal solution of $(S)_{\lam
^{**}+\delta}$ with $\delta >0$, then we would have for $\lam \le
\lam ^{**}$,
$$
-\Delta w=\frac{(\lam ^{**}+\delta )f(x)}{(1-w)^{2 }}\ge \frac{\lam
f(x)}{(1-w)^{2 }} \quad x\in \Omega \,.
$$
Since the minimal solutions $u_\lambda$ satisfy $ -\Delta u_{_\lam}=
\frac{\lam f(x)}{(1-u_{_\lam})^{2 }} \quad x\in \Omega \, $ for all
$0<\lam <\lam ^{**}$, it follows from Lemma 3.4 that $1> w\ge
u_{_\lam}$ for all $0<\lam <\lam ^{**}$. Consequently, $\b{u}=\lim
_{\lam \nearrow \lam ^{**}}u_{_\lam} $ would exist. Now from the
definition of $\lam ^{**}$ and Lemma 3.4, we must have $w\equiv
\b{u}$ and $\delta =0 $ on $\Omega $ which is  a contradiction, and
hence $\lam ^{**}= \lam ^*$.

(2) From the first part we conclude that if $0<\lam <\lam ^*$ and
$u=u_{_\lam}$, then the smallest eigenvalue of $-\Delta -\frac{2\lam
f(x)}{(1-u)^{3}}$ is positive. Applying the maximum principle, it is
easy to show that $u_{_\lam}(x)$ is increasing with respect to $\lam
$ (More details can be found in the proof of Theorem 1.2(1) below).
  That $\mu _{_{\lam }}$ is decreasing with respect to $\lam $ follows
now easily from the variational characterization of $\mu _{_{\lam
}}$ and the convexity of $(1-u)^{-3}$ with respect to $u$.
\endproof

\begin{rem} For the case where $f(x)>0$ on $\Omega $, Lemma 3 of
\cite{CR} gives $\mu (1,0)$ as an upper bound for $\lam ^{**}(=\lam
^{*})$. It is worth noting that our upper bound $\bar {\lam}$  in
Theorem 1.1 gives a better estimate. Indeed,if $f\equiv 1$ then $\mu
(1,0)=\mu_\Omega/2$ while the estimate in Theorem 1.1 gives
$\frac{4\mu_\Omega}{27}$ for an upper bound.
\end{rem}

\noindent  {\bf Proof of Theorem 1.2(1):} By Theorem 3.1, it
suffices
   to prove that for each $x\in \Omega $, the function $\lambda
\to u_{\lam}(x)$ is differentiable and strictly increasing on $(0,
\lambda^*)$. Setting $ F(\lam , u_\lambda (x))= -\Delta u_\lambda
-\frac{  \lam f(x)}{(1-u_\lambda)^{ 2}}\,,$ Proposition 3.3 then
implies that $F_{u_{\lam}}(\lam ,u_{_\lam})$ on $\Omega $ is
invertible for $0<\lam<\lam ^*$. It then follows from the Implicit
Function Theorem that $u_{_\lam}(x)$ is differentiable with respect
to $\lam$.

Consider now for  $\lam _1<\lam _2<\lam ^*$,  their corresponding
minimal positive solutions $u_{\lam_1}$ and $u_{\lam_2}$ and let
$u^*$ be a positive solution for $(S)_{\lam _2}$. For the monotone
increasing series $\{u_n(\lam _1; x)\}$ defined in $(\ref{3:1})$, we
then have $u^*>u_0(\lam _1; x)\equiv 0$, and if $u_{n-1}(\lam _1;
x)\le u^*$ in $\Omega $, then
$$\begin {array}{lll}
      -\Delta (u^*-u_n) = f(x)\big [\displaystyle\frac{\lam
_2}{(1-u^*)^{2}}- \displaystyle
      \frac{\lam _1}{(1-u_{n-1})^{2}}\big ]\ge 0
      \,, \quad x \in \Omega \,\\[3mm]
      \quad \ \ \ \ u^*-u_n=0 \,,\ \ \ \ \quad x \in \partial\Omega \,.
\end{array} $$
So we have $u_{n}(\lam _1; x)\le u^*$ in $\Omega $. Therefore,
$u_{\lam_1}=\lim _{n\to \infty}u_n(\lam _1; x)\le u^*$ in $\Omega $,
and in particular $u_{\lam_1}\le u_{\lam_2}$ in $\Omega $.
Therefore, $\frac{du_\lambda (x)}{d\lambda}\ge 0$ for all  $x\in
\Omega $.

Finally, by differentiating $(S)_{\lam }$ with respect to $\lam$ we
get
$$\begin {array}{ll}
\displaystyle-\Delta \frac{du_\lambda}{d\lambda} -
\displaystyle\frac{2 \lam
f(x)}{(1-u_{_\lam})^{3}}\frac{du_\lambda}{d\lambda} =
\displaystyle\frac{
f(x)}{(1-u_{_\lam})^{2 }}\ge 0 \,, \quad x \in \Omega \,\\[3mm]
\quad \quad \quad \quad \quad \quad  \quad \quad \  \,  \quad \quad
\ \ \ \ \displaystyle \frac{du_\lambda}{d\lambda} \ge 0 \,,\   \quad
x \in
\partial\Omega \,.
\end{array} $$
Applying the strong maximum principle, we conclude that
$\frac{du_\lambda}{d\lambda}
>0$ on $\Omega $ for all $0<\lam<\lam ^*$.
\endproof

\subsection{Energy estimates and regularity}

We start with the following easy observation.

\begin{lem}
Any positive (weak) solution $u$ in $ H^1_0(\Omega)$ of $(S)_{\lam}$
satisfies $\int _{\Omega }\frac{f}{(1-u)^2}dx<\infty $.
\end{lem}

\noindent  {\bf Proof:} Since $u\in H^1_0(\Omega )$ is a positive
solution of $(S)_{\lam}$, we have
$$
\int _{\Omega }\frac{f}{(1-u)^2}-\int _{\Omega }\frac{f}{1-u}=\int
_{\Omega }\frac{uf}{(1-u)^2}=\int _{\Omega }|\nabla u|^2<C\,,
$$
which implies that
$$
\int _{\Omega }\frac{f}{(1-u)^2}\le C+\int _{\Omega
}\frac{f}{1-u}\le C+\int _{\Omega }\big[C\varepsilon
\frac{f}{(1-u)^2}+\frac{C}{\varepsilon }f\big]\le C+C\varepsilon
\int _{\Omega }\frac{f}{(1-u)^2}
$$
with $\varepsilon >0$. Therefore, by choosing $\varepsilon >0$ small
enough, we conclude that $\int _{\Omega
}\frac{f}{(1-u)^2}<\infty\,.$
\endproof

That  $f/(1-u)\in L^2(\Omega )$ is unfortunately not sufficient to
obtain regularity results for the solutions. However, we now show
that the situation is much better if $f/(1-u)$ has better
integrability properties.

\begin{thm} For any bounded domain $\Omega \subset \R^N$ and any constant
$C>0$ there exists $0<K(C, N)<1$ such that a positive weak solution
$u$ of $(S )_\lambda$ $(0<\lambda <\lambda^*)$ is a classical
solution and $ \parallel u\parallel _{_{C(\Omega)}}\leq K(C, N)$
provided one of the following conditions holds:
\begin{enumerate}
\item  $N=1$ and   $\|\frac{f}{(1-u)^3}\|_{_{L^1(\Omega )}} \leq C$.
\item  $N=2$ and   $\|\frac{f}{(1-u)^3}\|_{_{L^{1+\epsilon}(\Omega )}}
\leq C$ for some $\epsilon>0$.
\item $N> 2$ and   $\|\frac{f}{(1-u)^3}\|_{_{L^{N/2}(\Omega )}} \leq C$.

\end{enumerate}
\end{thm}

\noindent  {\bf Proof:} We prove this lemma by considering the
following three cases separately:\\

\noindent (1)  If $N=1$, then for any $I>0$ we write using the
Sobolev inequality with constant $K(1)>0$,
\begin{equation}\arraycolsep=1.5pt \begin {array}{lll} \arraycolsep=1.5pt
   K(1)\parallel (1-u)^{-1}-1\parallel ^2_{L^{\infty} }&\le &
\displaystyle \int _{\Omega }\big |\nabla
[(1-u)^{-1}-1]\big |^2\\[4mm]
&=&\displaystyle\frac{1}{3}\displaystyle \int _{\Omega }\nabla
u\cdot \nabla \big[(1-u)^{-3}-1\big] \\[4mm]
&= &\displaystyle\frac{\lam}{3} \displaystyle \int _{\Omega
}f(1-u)^{-2}[(1-u)^{-3}-1]\\[4mm]
&\le  &  CI+C\displaystyle \int _{ \{(1-u)^{-3 } \ge I
\}}f(1-u)^{-5}\\[4mm]
& \le  &  CI+C\displaystyle \int _{ \{(1-u)^{-3 } \ge I
\}}8f(1-u)^{-2}\\[4mm]
&&  +C\displaystyle \int _{ \{(1-u)^{-3 } \ge I
\}}f\big[(1-u)^{-3}+2(1-u)^{-2}+4(1-u)^{-1}\big]\big[(1-u)^{-1}-1\big]^2\\[4mm]
   &\le   &  CI+C+C\parallel (1-u)^{-1}-1\parallel
^2_{L^{\infty}(\{(1-u)^{-3} \ge I \}) }\displaystyle \int _{
\{(1-u)^{-3} \ge I
\}}\frac{f}{(1-u)^3} \\[4mm]
&\le  &  CI+C+C\varepsilon (I)\parallel (1-u)^{-1}-1\parallel
^2_{L^{\infty} }\,
\end{array}
\label{4:93}
\end{equation}
with $ \varepsilon (I)=\displaystyle \int _{ \{(1-u)^{-3} \ge I
\}}\frac{f}{(1-u)^3}\,. $
  From the assumption $f/(1-u)^3\in L^1(\Omega )$, we have
$\varepsilon (I)\to 0$ as $I\to \infty$. We now choose $I$  such
that $\varepsilon (I)\leq \frac{K(1)}{2C}$, so that the above
estimates imply that $
   \parallel (1-u)^{-1}-1\parallel _{L^{\infty} }<K(C)\,.
$ Standard regularity theory for elliptic problems now imply that
$1/(1-u)\in  C^{2,\alp }(\Omega )$. Therefore, $u$ is classical and
there exists a constant $K(C, N)$ which can be taken strictly less
than $1$ such that $
\parallel u\parallel _{_{C(\Omega)}}\le K(C,N)<1$.\\

\noindent (2) The case when $N=2$ is similar as one can use that
$H^1_0$ embeds in $L^p$ for any $p<+\infty$.\\

\noindent (3) The case when $N>2$ is more elaborate and we first
show that $(1-u)^{-1}\in L^q(\Omega )$ for all $q\in (1,\infty )$.
Since $u\in H^1_0(\Omega )$ is a solution of $(S)_{\lam}$, we
already have
   $\int _{\Omega
}\frac{f}{(1-u)^2}<C.$  Now we proceed by iteration to show that if
$\int _{\Omega }\frac{f}{(1-u)^{2+2\theta }}<C$ for some $\theta
\geq 0$, then $\int _{\Omega }\frac{1}{(1-u)^{2^*(1+\theta )}}<C.$

Indeed, for any constant $\theta \ge 0$ and $\ell
>0$ we choose a test function $\phi =[(1-u)^{-3}-1]\min
\{(1-u)^{-2\theta }, \ell ^2\}$. By applying this test function to
both sides of $(S)_{\lam}$, we have
   \begin{equation}\arraycolsep=1.5pt \begin {array}{lll}
\arraycolsep=1.5pt &&\lam \displaystyle \int _{\Omega
}f(1-u)^{-2}[(1-u)^{-3}-1]\min \{(1-u)^{-2\theta }, \ell ^2\}=
\displaystyle \int _{\Omega }\nabla u\cdot \nabla
\big[\big((1-u)^{-3}-1\big)\min
\{(1-u)^{-2\theta }, \ell ^2\}\big]\\[4mm]
&&=3\displaystyle \int _{\Omega }|\nabla u|^2(1-u)^{-4}\min
\{(1-u)^{-2\theta }, \ell ^2\}+2\theta \displaystyle \int _{
\{(1-u)^{-\theta} \le \ell \}}|\nabla u|^2(1-u)^{-2\theta
-1}[(1-u)^{-3}-1] \,. \
\end{array}
\label{4:91}
\end{equation}
We now suppose $\int _{\Omega }\frac{f}{(1-u)^{2+2\theta }}<C.$ We
then obtain from $(\ref{4:91})$ and the fact that
$\frac{1}{(1-u)^5}\leq C_{_I}\frac{1}{(1-u)^3}(\frac{1}{1-u}-1)^2$
for $(1-u)^{-3} \geq I >1$ that
\begin{equation}\arraycolsep=1.5pt \begin {array}{lll} \arraycolsep=1.5pt
&&\displaystyle \int _{\Omega }\big |\nabla
[\big((1-u)^{-1}-1\big)\min
\{(1-u)^{-\theta }, \ell \}]\big |^2\\[4mm]
&&\leq \displaystyle 2\int _{\Omega }|\nabla u|^2(1-u)^{-4}\min
\{(1-u)^{-2\theta }, \ell ^2\}+2\theta ^2\displaystyle \int _{
\{(1-u)^{-\theta} \le
\ell \}}|\nabla u|^2(1-u)^{-2\theta -2}\big[(1-u)^{-1}-1\big]^2\\[4mm]
&&=\displaystyle 2\int _{\Omega }|\nabla u|^2(1-u)^{-4}\min
\{(1-u)^{-2\theta }, \ell
^2\}\\[4mm]
   &&\quad +2\theta ^2\displaystyle \int _{ \{(1-u)^{-\theta}
\le \ell \}}|\nabla u|^2(1-u)^{-2\theta -1}\big[(1-u)^{-3}-1+1+
(1-u)^{-1}-2(1-u)^{-2}\big]\\[4mm]
&&\le C\lam \displaystyle \int _{\Omega
}f(1-u)^{-2}[(1-u)^{-3}-1]\min
\{(1-u)^{-2\theta }, \ell ^2\}\\[4mm]
&&\le C\lam \displaystyle \int _{\Omega }f(1-u)^{-5}\min
\{(1-u)^{-2\theta }, \ell ^2\}\\[4mm]
&&\le CI+C\displaystyle \int _{ \{(1-u)^{-3 } \ge I
\}}f(1-u)^{-5}\min
\{(1-u)^{-2\theta }, \ell ^2\}\\[4mm]
&&\le CI+C\displaystyle \int _{ \{(1-u)^{-3 } \ge I
\}}f(1-u)^{-3}\big[(1-u)^{-1}-1\big]^2\min
\{(1-u)^{-2\theta }, \ell ^2\}\\[4mm]
&&\le CI+C\Big [\displaystyle \int _{ \{(1-u)^{-3} \ge I
\}}\Big (\frac{f}{(1-u)^3}\Big )^{\frac{N}{2}}\Big]^{\frac{2}{N}}\\[4mm]
&&\quad \times \Big [\displaystyle \int _{ \{(1-u)^{-3} \ge I
\}}\big ( \big[(1-u)^{-1}-1\big]\min \{(1-u)^{-\theta }, \ell
\}\big )^{\frac{2N}{N-2}}\Big]^{\frac{N-2}{N}}\\[4mm]
&&\le CI+C\varepsilon (I)\displaystyle \int _{\Omega }\big |\nabla
[\big((1-u)^{-1}-1\big)\min \{(1-u)^{-\theta }, \ell \}]\big |^2\,
\end{array}
\label{4:92}
\end{equation}
with
$$\varepsilon (I)=\Big [\displaystyle \int _{
\{(1-u)^{-3} \ge I \}}\Big (\frac{f}{(1-u)^3}\Big
)^{\frac{N}{2}}\Big]^{\frac{2}{N}}\,.
$$
  From the assumption $f/(1-u)^3\in L^{\frac{N}{2} }(\Omega )$ we have
$\varepsilon (I)\to 0$ as $I\to \infty$. We now choose $I$ such that
$\varepsilon (I)=\frac{1}{2C}$, and the above estimates imply that
$$
   \int _{ \{(1-u)^{-\theta } \le
\ell \}}\big |\nabla [(1-u)^{-\theta -1}-(1-u)^{-\theta }]\big
|^2\le CI\,,
$$
where the bound is uniform with respect to $\ell $. This estimate
leads to
$$
\arraycolsep=1.5pt \begin {array}{lll} \arraycolsep=1.5pt
\frac{1}{(\theta +1)^2}\int _{ \{(1-u)^{-\theta } \le \ell \}}\big
|\nabla [(1-u)^{-\theta -1}]\big |^2&=&\displaystyle\int _{
\{(1-u)^{-\theta } \le \ell \}}(1-u)^{-2\theta
-4}\big |\nabla u \big| ^2\\
&\le &CI+C\displaystyle\int _{ \{(1-u)^{-\theta } \le \ell
\}}(1-u)^{-2\theta -3}\big |\nabla u \big |
^2\\[4mm]
&\le &CI+\displaystyle\int _{ \{(1-u)^{-\theta } \le \ell \}}\big
[C\varepsilon (1-u)^{-2\theta -4}+C/\varepsilon \big
]\big |\nabla u \big | ^2\\[4mm]
&\le &CI+C\varepsilon \displaystyle\int _{ \{(1-u)^{-\theta } \le
\ell \}}(1-u)^{-2\theta -4}\big |\nabla u \big| ^2\,
\end{array}
$$
with $\varepsilon >0$. This means that for $\varepsilon >0$
sufficiently small
$$
\int _{ \{(1-u)^{-\theta } \le \ell \}}\big |\nabla (1-u)^{-\theta
-1}\big| ^2=\int _{ \{(1-u)^{-\theta } \le \ell \}}(\theta
+1)^2(1-u)^{-2\theta -4}\big |\nabla u \big| ^2<C\,.
$$
So we can let $\ell \to\infty $ and we get that $ (1-u)^{-\theta
-1}\in H^1(\Omega ) \hookrightarrow L^{2^*}(\Omega )$, which means
that $\int _{\Omega }\frac{1}{(1-u)^{2^*(1+\theta )}}<C.$

By iterating the above argument for  $\theta
_i+1=\frac{N}{N-2}(\theta _{i-1}+1)$ for $i\ge 1$ and starting with
$\theta _0=0$, we find that $1/(1-u)\in L^q(\Omega )$ for all $q\in
(1,\infty )$.

Standard regularity theory for elliptic problems applies again to
give that $1/(1-u)\in  C^{2,\alp }(\Omega )$. Therefore, $u$ is a
classical solution and there exists a constant $0<K(C,N)<1$ such
that $
\parallel u\parallel _{_{C(\Omega)}}\le K(C,N)<1$. This completes the
proof of Theorem 3.6.
\endproof

\begin{thm} For any dimension $1\leq N<8$, there exists a constant
$0<C(N)<1$ independent of $\lambda$ such that for any $0<\lam <\lam
^*$, the  minimal solution $u_\lambda$  satisfies $\parallel
u_{_\lam} \parallel _{_{C(\Omega)}}\le C(N)$.
\end{thm}

The theorem, which gives Theorem 1.2(2), will follow from the
following uniform energy estimate on the minimal solutions
$u_{_\lam}$.

\begin{lem} There exists a constant $C(p)>0$ such that for each
$\lam\in (0,\lam ^*)$, the minimal solution $u_{_\lam}$ satisfies
$\|\frac{f}{(1-u_{_\lam})^3}\|_{L^{p}(\Omega)}\leq C(p)$ as long as
$p < 1+\frac{4}{3}+2\sqrt{\frac{2}{3}}$.
\end{lem}

\noindent  {\bf Proof:} Proposition 3.3 implies that
\begin{equation}
\lam \int _{\Omega }\frac{2f(x)}{(1-u_{_\lam})^{3}}w^2dx\le -\int
_{\Omega }w\Delta wdx=\int _{\Omega }|\nabla w|^2dx\,, \label{4:1}
\end{equation}
for all $0<\lam<\lam ^*$ and nonnegative $w\in H^1_0(\bar
{\Omega})$. Setting
\begin{equation}
w=(1-u_{_\lam})^i-1>0\,, \  where \quad -2 -\sqrt{6 }<i<0\,,
\label{4:2}
\end{equation}
then $(\ref{4:1})$ becomes
\begin{equation}
i^2\int _{\Omega }(1-u_{_\lam})^{2i-2}|\nabla u_{_\lam}|^2dx\ge \lam
\int _{\Omega }\frac{2
[1-(1-u_{_\lam})^i]^2f(x)}{(1-u_{_\lam})^{3}}dx\,. \label{4:3}
\end{equation}
On the other hand, multiplying $(S)_{\lam}$ by
$\frac{i^2}{1-2i}[(1-u_{_\lam})^{2i-1}-1]$ and applying integration
by parts yield that
\begin{equation}
i^2\int _{\Omega }(1-u_{_\lam})^{2i-2}|\nabla u_{_\lam}|^2dx =\lam
\frac{i^2}{2i-1}\int _{\Omega }
\frac{[1-(1-u_{_\lam})^{2i-1}]f(x)}{(1-u_{_\lam})^{2}}dx\,.
\label{4:4}
\end{equation}
And hence $(\ref{4:3})$ and $(\ref{4:4})$ reduce to
\begin{equation}\arraycolsep=1.5pt \begin {array}{lll} \arraycolsep=1.5pt
&&\displaystyle \frac{\lam \, i^2}{2i-1}\displaystyle \int _{\Omega
} \frac{f(x)}{(1-u_{_\lam})^{2}}dx-2\lam\displaystyle \int _{\Omega
} \frac{f(x)}{(1-u_{_\lam})^{3}}dx+4\lam  \displaystyle \int
_{\Omega }
\frac{f(x)}{(1-u_{_\lam})^{3-i}}dx\\[3mm]
&&\ge \lam (2+\displaystyle \frac{i^2}{2i-1})\displaystyle \int
_{\Omega } \frac{f(x)}{(1-u_{_\lam})^{3-2i}}dx\,.
\end{array} \label{4:5}\end{equation}
  From the choice of $i$ in $(\ref{4:2})$ we have $2+\frac{
i^2}{2i-1}>0$. So $(\ref{4:5})$ implies that
\begin{equation}
\arraycolsep=1.5pt \begin {array}{lll} \displaystyle \int _{\Omega
}\frac{f(x)}{(1-u_{_\lam})^{3-2i}}dx&\le &C\displaystyle
\int _{\Omega }\frac{f(x)}{(1-u_{_\lam})^{3-i}}dx\\[3mm]
&\le &C\Big(\displaystyle \int _{\Omega }\Big |
\frac{f^{\frac{3-i}{3-2i}}}{(1-u_{_\lam})^{3-i}}\Big
|^{\frac{3-2i}{3-i}} dx\Big)^{\frac{3-i}{3-2i}}\cdot
\Big(\displaystyle \int _{\Omega }\Big |
    f^{\frac{-i}{3-2i}}\Big
|^{\frac{3-2i}{-i}} dx\Big)^{\frac{-i}{3-2i}}\\[3mm]
&\le &C\Big (\displaystyle \int _{\Omega
}\frac{f(x)}{(1-u_{_\lam})^{3-2i}}dx\Big)^{\frac{3-i}{3-2i}}\,,
\end{array}
\label{4:6}\end{equation} where Holder's inequality is applied. From
the above we deduce that
\begin{equation}
\int _{\Omega }\frac{f(x)}{(1-u_{_\lam})^{3-2i}}dx\le C\,.
\label{4:7}\end{equation} Further we have
\begin{equation}
\arraycolsep=1.5pt \begin {array}{lll} \arraycolsep=1.5pt
\displaystyle \int _{\Omega }\Big |
\frac{f(x)}{(1-u_{_\lam})^{3}}\Big |^{\frac{3-2i}{3}}
dx&=&\displaystyle \int _{\Omega }f^{\frac{-2i}{3}}\cdot
\frac{f}{(1-u_{_\lam})^{3-2i}}dx\\[3mm]
&\le &C\displaystyle \int _{\Omega
}\frac{f}{(1-u_{_\lam})^{3-2i}}dx\le C\,. \
\end{array}
\label{4:8}\end{equation} Therefore, we get that \bsub \label{4:9}
\begin{equation}
\parallel \frac{f(x)}{(1-u_{_\lam})^{3}}\parallel
_{_{L^p}}\le C \,, \label{4:9a}
\end{equation}
where --in view of  $(\ref{4:2})$--
\begin{equation}
       p=\frac{3-2i}{3}\le 1+\frac{4}{3}+2\sqrt{\frac{2}{3}}\,.
\label{4:9b}
\end{equation}
\esub

\endproof

\noindent {\bf Proof of Theorem 3.7:} This follows from Lemma 3.8
and Theorem 3.6, where $p=1$ when the dimension $N=1$, $p$ can be
taken to be $1+\frac{4}{3}$ when $N=2$. For $N>2$, the reasoning
applies as long as $\frac{N}{2} < 1+\frac{4}{3}+2\sqrt{\frac{2}{3}}$
which happens when $N<8$.\endproof

Finally, we note the following easy comparison results and we omit
the details.
\begin{cor}
Suppose  $f_1, f_2:\Omega \to (0,1]$ satisfy  $f_1(x)\leq f_2(x)$ on
$\Omega $, then $\lambda^*(\Omega, f_1) \ge  \lambda^*(\Omega, f_2)$
and for $0<\lam < \lambda^*(\Omega, f_1)$ we have  $u_1(\lam, x
)\leq u_2(\lam, x )$ on $\Omega$, where $u_1(\lam ,x)$ (resp.,
$u_2(\lam ,x)$) are the unique minimal positive solution of
$$
\hbox{ $-\Delta u = \frac{\lambda f_1(x)}{(1-u)^2}$ (resp.,
     $
-\Delta u = \frac{\lambda f_2(x)}{(1-u)^2})$ on $\Omega$ and  $u =
0$  on $\partial\Omega. $}
    $$
Moreover, if $f_1(x)> f_2(x)$ on a subset of positive measure, then
$u_1(\lam ,x)<u_2(\lam ,x)$ for all $x\in \Omega $.
\end{cor}

We note that if one considers the cases of power-law or exponential
profiles for $(S)_{\lam }$ defined in a ball, then the minimal
positive solution corresponds to the lowest branch in the
bifurcation diagram, the one connecting the origin point $\lam =0$
to the first fold at $\lam =\lam ^*$, see section \S 5.

\subsection{Existence of solutions at $\lam =\lam ^*$}

In this subsection, we study the existence of positive solutions at
the critical voltage $\lam =\lam ^*$. We first deal with the
existence of minimal solutions $u_{\lam ^*}$ for $(S)_{\lam ^*}$.

\begin{lem}
Suppose there exists $0<C<1$ such that $\parallel u_{_\lam}\parallel
_{_{C(\bar {\Omega})}}\le C$ for each $\lam <\lam ^*$. Then
    $u_{_{\lam^*}}=\lim _{_{\lam \nearrow \lam ^*}}u_{_\lam}$ exists in
the $C^{2,\alp}(\bar {\Omega})$ topology for some $0<\alp <1$.
Moreover, there exists $\delta >0$ such that the solutions of
$(S)_{\lam}$ near $(\lam ^*, u_{_{\lam^*}})$ form a curve $\rho (s)
=\{(\bar\lam (s),v(s)):\, |s|<\delta\}$, and the pair $(\bar\lam
(s),v(s))$ satisfies:
\begin{equation}\bar\lam (0)=\lam ^*,\ \bar\lam '(0)=0,\ \bar\lam
''(0)<0\,, \ \ {\rm and} \ \ v(0)=u_{_{\lam^*}},\ v'(0)(x)> 0 \ {\rm
in}\ \Omega \,.\label{4:95}
\end{equation}
\end{lem}

\noindent  {\bf Proof:}  The proof is similar to  a related result
of Crandall and Rabinowitz cf.~\cite{CR1} \cite{CR}, so we will be
brief. Firstly, the assumed upper bound on $u_{_\lam}$ in $C^1$ and
standard regularity theory, show that if $f\in C(\bar \Omega)$ then
$\parallel u_{_\lam}\parallel _{_{C^{2,\alp }(\bar {\Omega})}}\le C$
for some $0<\alp <1$ (while if $f\in L^\infty$, then  $\parallel
u_{_\lam}\parallel _{_{C^{1,\alp }(\bar {\Omega})}}\le C$).  It
follows that $\{(\lam ,u_{_\lam})\}$ is precompact in the space
$\mathbb{R}\times C^{2,\alp }$, and hence we have a limiting point
$(\lam ^*,u_{_{\lam^*}})$ as desired. Since $\frac{\lam ^*
f(x)}{(1-u_{_{\lam^*}})^2}$ is nonnegative, Theorem 3.2 of
\cite{CR1} characterizes the solution set of $(S)_{\lam}$ near
$(\lam ^*,u_{_{\lam^*}})$: $\bar\lam (0)=\lam ^*$, $\bar\lam
'(0)=0$, $v(0)=u_{_{\lam^*}}$ and $v'(0)>0$ in $\Omega$. Finally,
the same computation as in Theorem 4.8 in \cite{CR1} gives that
$\bar\lam ''(0)<0$.
\endproof

\begin{rem}
Lemma 3.10 implies that if the minimal solution $u_{_\lam}(x)$
satisfies $\parallel u_{_\lam}\parallel _{_{C(\bar {\Omega})}}\le
C<1$ (which occurs when $N<8$),  then there exists two distinct
solutions for $(S)_{\lam}$ for $\lam$ in a deleted left neighborhood
of $\lam ^*$. A version of this result will be established
variationally in an upcoming paper.
\end{rem}

The following theorem gives the uniqueness of (classic) solutions
for $(S)_{\lam ^*}$.


\begin{thm}Suppose
there exists $0<C<1$ such that $\parallel u_{_\lam}\parallel
_{_{C(\bar {\Omega})}}\le C$ for each $\lam <\lam ^*$. Then the
minimal solution $u_{_{\lam^*}}=\lim _{_{\lam \nearrow \lam
^*}}u_{_\lam}$ obtained above satisfies the following properties:
\begin{enumerate}
\item The smallest eigenvalue $\mu (\lam )$ at $\lam =\lam ^*$ of the
linearized operator $L_{u_{_{\lam }}, \lam }=-\Delta -\frac{2\lam
f(x)}{(1-u_{_{\lam }})^3}$ on $\Omega $ is zero.
\item $u_{_{\lam
^*}}$ is the unique solution of $(S)_{\lam ^*}$.
\end{enumerate}
\end{thm}

\noindent  {\bf Proof:} (1) Applying Proposition 3.3(2) we see that
$\mu (\lam )>0$ on the minimal branch for any $\lam <\lam ^*$, hence
the limit $\mu (\lam ^*)\ge 0$. If now $\mu (\lam ^*)>0$ the
Implicit Function Theorem could be applied to the operator
$L_{u_{_{\lam ^*}}, \lam ^*}$, and would allow the continuation of
the minimal
   branch $\lam \mapsto u_{_{\lam }}$ of classical solutions beyond
$\lam ^*$, which is a contradiction and hence $\mu (\lam ^*)=0$.

(2)\  Suppose now $u$ is any solution such that  $u\ge u_{_{\lam
^*}}$. Since $\mu (\lam ^*)=0$, let $\phi$ be any
    positive eigenfunction  in the kernel of
$L_{u_{\lambda^*},\lambda^*}$ and write,
$$
-\phi \Delta (u-u_{_{\lam ^*}}) =\lam
^*f(x)\big[\frac{1}{(1-u)^2}-\frac{1}{(1-u_{_{\lam ^*}})^2}\big]\phi
\,,
$$
which yields that
$$
-\int _{\Omega }(u-u_{_{\lam ^*}})\Delta \phi =\lam ^*\int _{\Omega
}f(x)\big[\frac{1}{(1-u)^2}-\frac{1}{(1-u_{_{\lam
^*}})^2}\big]\phi\,.
$$
On the other hand, since $ -\Delta \phi =\frac{2\lam
^*f(x)}{(1-u_{_{\lam ^*}})^3}\phi $, we have
$$
\lam ^*\int _{\Omega
}f(x)\big[\frac{1}{(1-u)^2}-\frac{1}{(1-u_{_{\lam
^*}})^2}-\frac{2}{(1-u_{_{\lam ^*}})^3}(u-u_{_{\lam ^*}})\big]\phi
=0\,.
$$
Since the integrand is nonnegative it follows that
\begin{equation}
\frac{1}{(1-u)^2}=\frac{1}{(1-u_{_{\lam
^*}})^2}+\frac{2}{(1-u_{_{\lam ^*}})^3}(u-u_{_{\lam ^*}})\quad a.e.\
in\ \Omega\,. \label{4:96}
\end{equation}
If now $\parallel u\parallel _{_{L^{\infty }}}\le C<1$, then $u$ is
a classical solution as in Theorem 3.6, and we conclude that
$u\equiv u_{_{\lam ^*}}$ on $\Omega $.
\endproof

Now Theorem 1.3 is a direct result of Theorem 1.2(2) (or Theorem
3.7), Lemma 3.10 and Theorem 3.11.

\section{Uniqueness and Multiplicity of Solutions}
The purpose of this section is to discuss uniqueness and
multiplicity of solutions for $(S)_{\lam}$. Note that Lemma 3.10
gives that for some $0<\lam ^*_2<\lam ^*$, there exists at least two
solutions for $(S)_{\lam}$ with $\lam \in (\lam ^*_2,\lam ^*)$,
which is Theorem 1.4(2). In the following we shall focus on the
uniqueness when $\lam$ is small enough. We first define non-minimal
solutions for $(S)_{\lam}$ as follows:

\begin{defn} A solution $0\le u<1$ is said to be a non-minimal positive
solution of $(S)_{ \lam}$, if there exists another positive solution
$v$ of $(S)_{ \lam}$ and a point  $x\in \Omega$ such that
$u(x)>v(x)$.
\end{defn}

\begin{lem} Suppose $u$ is a non-minimal solution of $(S)_\lambda$
with $\lam \in (0,\lam^*)$. Then the smallest eigenvalue $\mu (\lam
)$ of the linearized operator $L_{u, \lam }=-\Delta -\frac{2\lam
f(x)}{(1-u)^3}$ on $\Omega $ must be negative.
\end{lem}

\noindent  {\bf Proof:} For any fixed $\lam \in (0,\lam^*)$, let
$u_{_\lam}$ be the minimal solution of $(S)_\lam$. Clearly we have
$w=u-u_{_\lam}\ge 0$ in $\Omega $, and
$$
-\Delta w-\frac{\lam (2-u-u_{_\lam})f}{(1-u)^2(1-u_{_\lam})^2}w= 0\,
\quad {\rm in}\ \Omega\,.
$$
Hence we deduce from the strong maximum principle that $u_{_\lam}<
u$ in $\Omega $.

Let $\Omega _0=\{x\in \Omega: f(x)=0\}$ and $\Omega /\Omega
_0=\{x\in \Omega: f(x)>0\}$. Direct calculations give that
\begin{equation}
-\Delta (u-u_{_\lam})-\frac{2\lam f}{(1-u)^3}(u-u_{_\lam})=\lam
f\big[\frac{1}{(1-u)^2}-\frac{1}{(1-u_{_\lam
})^2}-\frac{2}{(1-u)^3}(u-u_{_\lam})\big]=\begin{cases}
   0\,,\, & x\in \Omega _0\,;\\[2mm]
<0\,,\, & x\in \Omega /\Omega _0\,.  \end{cases} \label{4:13}
\end{equation}
  From this we get
\begin{equation}\lam
\displaystyle\int _{\Omega /\Omega
_0}f\big[\displaystyle\frac{1}{(1-u)^2}-\displaystyle\frac{1}{(1-u_{_\lam
})^2}-\displaystyle\frac{2}{(1-u)^3}(u-u_{_\lam})\big](u-u_{_\lam})<0\,.
\label{4:11}
\end{equation}
Now suppose that $\mu (\lam )\ge 0$. Then for each $\phi \in
H^1_0(\Omega )$ we have
\begin{equation}
\langle L_{u, \lam }\phi ,\phi \rangle=\int _{\Omega }(|\nabla \phi
|^2 -\frac{2\lam f(x)}{(1-u)^3}\phi ^2)\ge 0\,. \label{4:12}
\end{equation}
Putting $\phi =u-u_{_\lam}$ in $(\ref{4:12})$, we get from the left
equality of $(\ref{4:13})$ that
$$\arraycolsep=1.5pt\begin{array}{lll}\arraycolsep=1.5pt&&\lam
\displaystyle\int _{\Omega /\Omega
_0}f\big[\displaystyle\frac{1}{(1-u)^2}-\displaystyle\frac{1}{(1-u_{_\lam
})^2}-\displaystyle\frac{2}{(1-u)^3}(u-u_{_\lam})\big](u-u_{_\lam})\\[4mm]
&&=\lam \displaystyle\int _{\Omega
}f\big[\displaystyle\frac{1}{(1-u)^2}-\displaystyle\frac{1}{(1-u_{_\lam
})^2}-\displaystyle\frac{2}{(1-u)^3}(u-u_{_\lam})\big](u-u_{_\lam})\ge
0\,
\end{array}$$
which contradicts $(\ref{4:11})$, and we are done.
\endproof

\begin{rem}  Proposition 3.3 and Lemma 4.1 give that for $0<\lam
<\lam ^*$, the smallest eigenvalue of $L_{u, \lambda} $ is
necessarily negative if $u$ is a non-minimal solution of
$(S)_{\lam}$, while the smallest eigenvalue of $L_{u, \lambda} $ is
positive if $u=u_{_\lam}$ is a minimal solution of $(S)_{\lam}$. For
parabolic problems of the type $(1.1)$, it is well-known that the
spectrum of the linearized operator about any steady-state solution
determines the stability of solutions for $(1.1)$. Therefore,
$u_{_\lam}(x)$ is the unique stable steady-state of $(1.1)$. In our
upcoming paper \cite{GG2}, we shall prove that the dynamic solution
of $(1.1)$ with $\lam <\lam^*$ (and $\lambda=\lambda^*$ for $N<8$)
will globally converge to its unique minimal solution
$u_{_\lam}(x)$.
\end{rem}

Now we are able to prove the following uniqueness result, which
completes the proof of Theorem 1.4.

\begin{thm}  For every $M>0$ there exists $0<\lam^* _1(M)<\lam ^*$
such that for $\lam \in (0,\lam^* _1(M))$ the equation $(S)_\lambda$
has a  unique  solution $v$ satisfying:
\begin{enumerate}
\item  $\|\frac{f}{(1-v)^{3}}\|_1 \leq M$ as long as the dimension $N=1$.
\item  $\|\frac{f}{(1-v)^{3}}\|_{1+\epsilon} \leq M$ and $N=2$.
\item  $\|\frac{f}{(1-v)^{3}}\|_{N/2} \leq M$ and $N>2$.

\end{enumerate}

\end{thm}
\noindent  {\bf Proof:} For any fixed $\lam \in (0,\lam^*)$, let
$u_{_\lam}$ be the minimal solution of $(S)_\lam$ and suppose
$(S)_{\lam}$ has a non-minimal solution $u$. Lemma 4.2 then gives
$$
\int _{\Omega}|\nabla (u-u_{_\lam})|^2dx<\int _{\Omega}\frac{2\lam
(u-u_{_\lam})^2f(x)}{(1-u)^3}dx\,.
$$
This implies in the case where $N>2$ that
$$\arraycolsep=1.5pt\begin{array}{lll}\arraycolsep=1.5pt
C(N)\Big (\displaystyle\int
_{\Omega}(u-u_{_\lam})^{\frac{2N}{N-2}}dx\Big)^{\frac{N-2}{N}}&<&
\lam
\displaystyle\int _{\Omega}\frac{2f(x)}{(1-u)^3} (u-u_{_\lam})^2dx\\[3mm]
&\le& 2\lam\Big( \displaystyle\int
_{\Omega}\big|\displaystyle\frac{f}{(1-u)^{3}}\big|^{\frac{N}{2}}
\Big )^{\frac{2}{N}} \Big(\int _{\Omega}
(u-u_{_\lam})^{\frac{2N}{N-2}}dx\Big)^{\frac{N-2}{N}}\\[4mm]
&\leq& 2\lam  M^{\frac{2}{N}}\Big(\displaystyle\int _{\Omega}
(u-u_{_\lam})^{\frac{2N}{N-2}}dx\Big)^{\frac{N-2}{N}}
\end{array}$$
which is a contradiction if $\lambda <\frac{C(N)}{2M^{\frac{2}{N}}}$
unless $u\equiv u_{_\lam}$. If $N=1$, then we write
$$
C(1)\|(u-u_{_\lam})\|_\infty^2 < \lam \int
_{\Omega}\frac{2f(x)}{(1-u)^3} (u-u_{_\lam})^2dx \le 2\lam
\|(u-u_{_\lam})\|^2_\infty \int _{\Omega}\frac{f}{(1-u)^3}dx
$$
   and the proof follows. A similar proof holds for dimension $N=2$.
\endproof

\begin{rem} The above gives uniqueness for small $\lambda$ among all
solutions that either stay away from $1$ or those that approach it
slowly. We do not know whether if $\lambda$ is small enough, any
positive solution $v$ of $(S)_\lambda$ satisfy $ \int
_{\Omega}(1-v)^{-\frac {3N}{2}}dx \leq M$ for some uniform bound $M$
independent of $\lambda$. Numerical computations do show that we may
have uniqueness for small $\lambda$ --at least for radially
symmetric solutions-- as long as $N\geq 2$.

\end{rem}

\section{Steady-State: Case of Power-Law Profile}

The issues of  uniqueness and multiplicity of solutions for $(S)_{
\lam }$ with $0<\lam <\lam ^*$, and even mere existence for $(S)_{
\lam ^*}$ with $N\ge 8$ seem to be quite challenging problems. In
this section, we discuss these problems in the case where $f$ has a
power-law permittivity profile, $i.e$, $f(x)=|x|^{\alp }\ (\alp \ge
0)$. We shall also consider the domain $\Omega$ to be a unit ball
$B_1(0)\subset \R ^N\ (N\ge 1)$ and $\lam \in (0, \lam ^*]$. In this
special case, the solutions of $(S)_{ \lam }$ must be radially
symmetric, and $(S)_{ \lam }$ is then reduced to the following
problem
\begin{equation}
\begin{array}{lll}
-u_{rr}-\displaystyle\frac{N-1}{r}u_r=\displaystyle\frac{\lam
r^{\alp
}}{(1-u)^2}\,,\quad 0<r\le 1,\\[3mm]
u'(0)=0\,,\quad u(1)=0\,.
\end{array}
     \label{1:8} \end{equation}
Here $r=|x|$ and $0<u=u(r)<1$ for $0<r<1$.

Looking first for a solution of the form
$$
u(r)=1-\beta w(P)\quad \mbox{with}\quad P=\gamma r\,,
$$
where $\gam$, $\beta >0$, equation $(\ref{1:8})$ implies that
$$
\gam ^2\beta \big (w''+\frac{N-1}{P}w'\big )=\frac{\lam P^{\alp
}}{\beta ^2\gam ^{\alp }}\frac{1}{w^2}\,.
$$
We set $w(0)=1$ and $\lam=\gam ^{2+\alp }\beta ^3$. This yields the
following initial value problem
\begin{equation}
\begin{array}{lll}
w''+\displaystyle\frac{N-1}{P}w'=\frac{P^{\alp }}{w^2}\,,\quad
P>0\,,\\[3mm]
w'(0)=0\,,\quad w(0)=1\,.
\end{array}
     \label{1:9} \end{equation}
Since $u(1)=0$ we have $\beta =1/w(\gam )$. Therefore, we conclude
that
\begin{equation}\left \{
\begin{array}{lll}
u(0)=1-\displaystyle\frac{1}{w(\gam )}\,,\\[3mm]
\lam=\displaystyle\frac{\gam ^{2+\alp }}{w^3(\gam )}\,,
\end{array}\right.
     \label{1:10} \end{equation}
where $w(\gam )$ is a solution of $(\ref{1:9})$.

As done in $\cite{P1}$,  one can numerically integrate the initial
value problem $(\ref{1:9})$ and use the results to compute the
complete bifurcation diagram for $(\ref{1:8})$. We show such a
computation of $u(0)$ versus $\lam$ defined in $(\ref{1:10})$ for
the slab domain $(N=1)$ in Fig.~2. In this case, one observes from
the numerical results that when $N=1$,
     \begin{itemize}
\item   There exists a unique solution for $(S)_{ \lam ^*}$;
\item For $0\le \alp \le 1$, there exist exactly two solutions for
$(S)_{ \lam }$ whenever $\lam \in (0,\lam ^*)$;
\item  For $\alp >1$, it is however difficult to see in any other
case the bifurcation diagram as $u(0)\to 1$.
\end{itemize}

This leads us to the question of determining  the asymptotic
behavior of $w(P)$ as $P\to \infty $. Towards this end, we proceed
it as follows.

\begin{figure}[htbp]
\begin{center}
{\includegraphics[width = 11cm,height=5.5cm,clip]{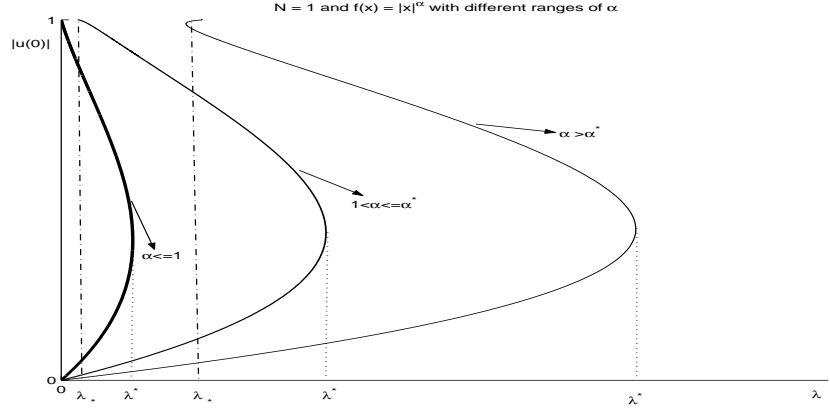}
\label{fig:fig2_2}}
  \caption{{\em   Plots of $|u(0)|$ versus $\lam$ for the power-law
permittivity profile $f(x)=|x|^{\alp }\ (\alp \ge 0)$ defined in the
slab domain ($N=1$). The numerical experiments show that it seems to
exist a constant $\alp ^*>1$ (analytically given in $(\ref{1:15})$),
such that the bifurcation diagrams are greatly different for
different ranges of $\alp$: $0\le \alp \le 1$, $1<\alp \le \alp ^*$
and $\alp >\alp ^*$.} } \label{fig:5:1}
\end{center}
\end{figure}

Set $\eta =logP$, $w(P)=P^BV(\eta )>0$ for some positive constant
$B$. Then we obtain from $(\ref{1:9})$ that
\begin{equation}
P^{B-2}V''+(2B+N-2)P^{B-2}V'+B(B+N-2)P^{B-2}V=\frac{P^{\alp
-2B}}{V^2}\,. \label{1:11} \end{equation} Choosing $B-2=\alp -2B$ so
that $B=(2+\alp )/3$, we get that
\begin{equation}
V''+\frac{3N+2\alp -2}{3}V'+\frac{(2+\alp )(3N+\alp
-4)}{9}V=\frac{1}{V^2}\,. \label{1:12} \end{equation} We notice from
$(\ref{1:12})$ that only for the case where $N\ge 2$ or $N=1$ with
$\alp >1$, that the equilibrium point $V_e$ of $(\ref{1:12})$ must
be positive and satisfies
\begin{equation}
V^3_e=\frac{9}{(2+\alp )(3N+\alp -4)}>0\,. \label{1:13}
\end{equation} When $N=1$, this is consistent with the numerical
observation of Fig.~2. Linearizing around this equilibrium point by
writing
$$
V=V_e+Ce^{\sig\eta }\,,\quad 0<C<<1,
$$
we obtain that
$$
\sig ^2+\frac{3N+2\alp -2}{3}\sig +\frac{(2+\alp)(3N+\alp
-4)}{3}=0\,.
$$
This reduces to \bsub \label{1:14}
\begin{equation}
\sig _{\pm }=-\frac{3N+2\alp -2}{6}\pm \frac{\sqrt {\bigtriangleup
}}{6}  \,, \label{1:14a}
\end{equation}
with
\begin{equation}
         \bigtriangleup =-8\alp ^2-(24N-16)\alp +(9N^2-84N+100) \,.
\label{1:14b}
\end{equation}
\esub We note that $\sig _{\pm }<0$ whenever $\bigtriangleup \ge 0$.
Define now
\begin{equation}
     \alp ^*=-\frac{1}{2}+\frac{1}{2}\sqrt{27/2} \,,\quad \alp
^{**}=\frac{4-6N+3\sqrt{6}(N-2)}{4} \ \
     (N\ge 8)\,.
\label{1:15}
\end{equation}
Next, we discuss on $N$ and $\alp $ by considering the sign of
$\bigtriangleup $:

\vspace {.1cm} {\bf Case (1).}\ \ $N$ and $\alp $ satisfy either one
of the followings: \bsub \label{1:16}
\begin{equation}
N=1  \quad with\quad 1<\alp \le \alp ^*\,; \label{1:16a}
\end{equation}
\begin{equation}
         N\ge 8 \quad with\quad 0\le \alp \le \alp ^{**}\,.  \label{1:16b}
\end{equation}
\esub

In this case, we have $\bigtriangleup \ge 0$ and
$$
V\sim \Big(\frac{9}{(2+\alp)(3N+\alp -4)}\Big)^{\frac{1}{3}}+
C_1e^{-\frac{3N+2\alp -2-\sqrt {\bigtriangleup }}{6}\eta }+\cdots
\,,\quad \mbox{as} \quad \eta \to +\infty \,.
$$
Further, we conclude that
$$
w\sim P^{\frac{2+\alp }{3}}\Big(\frac{9}{(2+\alp)(3N+\alp
-4)}\Big)^{\frac{1}{3}}+ C_1P^{-\frac{N-2}{2}+\frac{\sqrt
{\bigtriangleup }}{6}}+\cdots \,, \quad \mbox{as} \quad P \to
+\infty \,.
$$
Since $\lam =\gam ^{2+\alp }/w^3(\gam )$, we have
\begin{equation}
\lam \sim \lam _*=\frac{(2+\alp)(3N+\alp -4)}{9}\quad as \quad \gam
\to \infty \,.  \label{1:17}
\end{equation}
For this case, we compute the numerical results for $(\ref{1:16a})$
in Fig.~2 and for $(\ref{1:16b})$ in Fig.~3(b), respectively.

\vspace {.1cm} {\bf Case (2).}\ \ $N$ and $\alp $ satisfy any one of
the following three: \bsub \label{1:18}
\begin{equation}
N=1  \quad with\quad \alp > \alp ^*\,; \label{1:18a}
\end{equation}
\begin{equation}
        2\le N\le 7 \quad with\quad \alp \ge 0\,; \label{1:18b}
\end{equation}
\begin{equation}
       N\ge 8 \quad with \quad \alp > \alp ^{**}\,.  \label{1:18c}
\end{equation}
\esub

\begin{figure}[htbp]
\begin{center}
{\includegraphics[width = 8cm,height=5.5cm,clip]{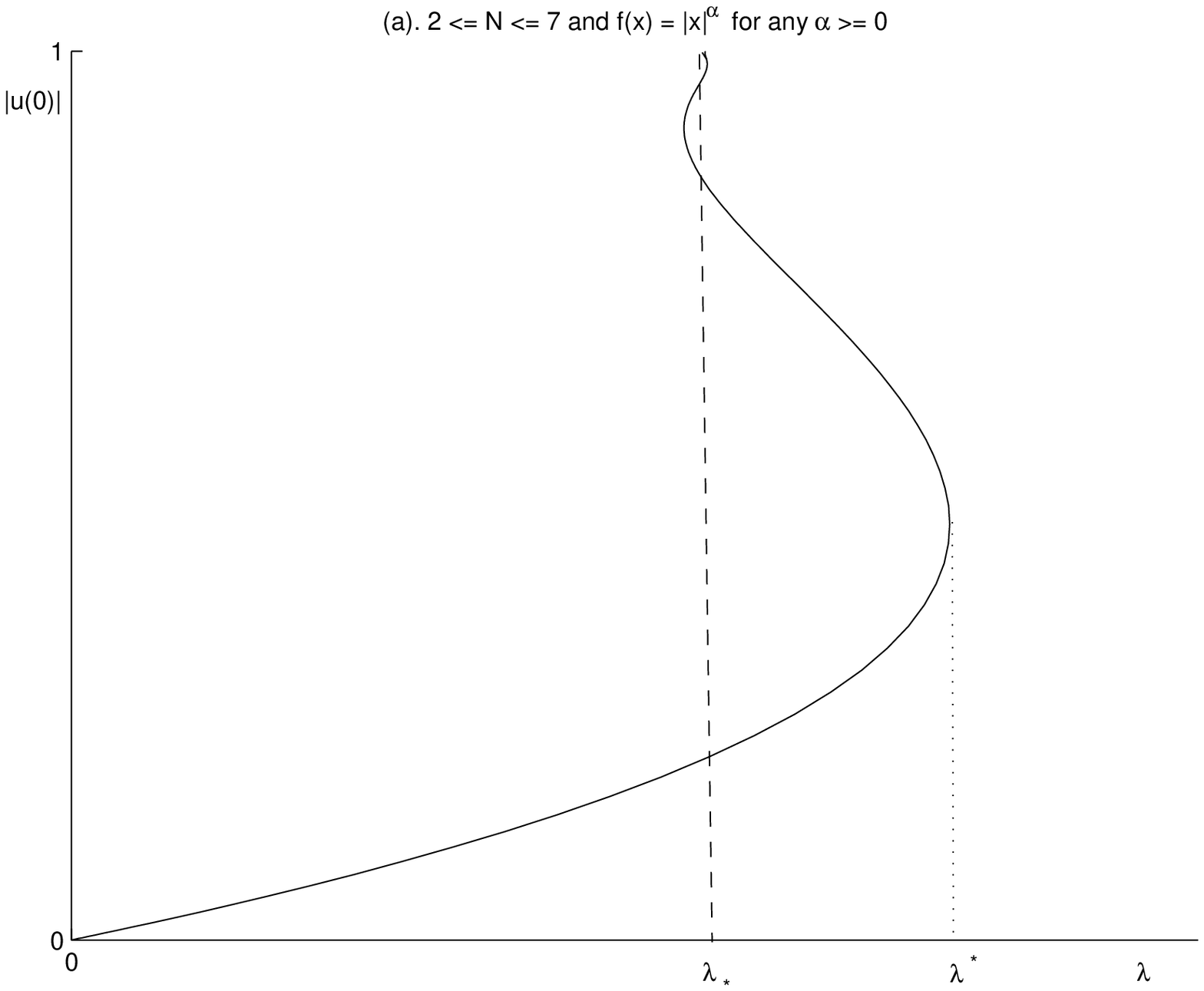}
\label{fig:fig2_2a}}
{\includegraphics[width = 8cm,height=5.5cm,clip]{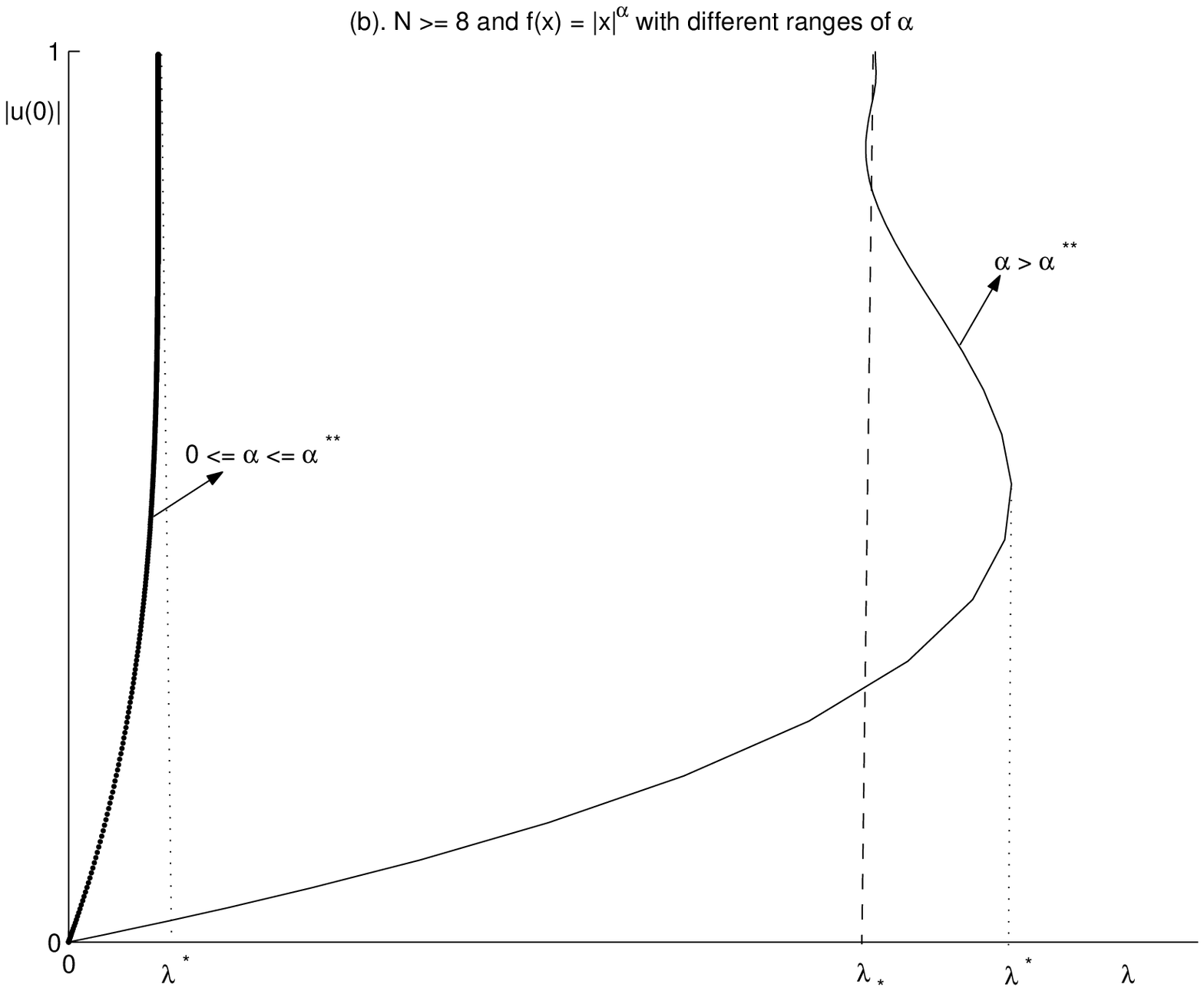}
\label{fig:fig2_2b}} \caption{{\em Left figure: Plots of $|u(0)|$
versus $\lam$ for the power-law permittivity profile $f(x)=|x|^{\alp
}\ (\alp \ge 0)$ defined in the unit ball $B_1(0)\subset \R^N$ with
$2\le N\le 7$. In this case, $|u(0)|$ oscillates around the value
$\lam _* $ defined in $(\ref{1:17})$ and there exists a unique
solution for $(S)_{ \lam ^*}$. Right figure: Plots of $|u(0)|$
versus $\lam $ for the power-law permittivity profile
$f(x)=|x|^{\alp }\ (\alp \ge 0)$ defined in the unit ball
$B_1(0)\subset \R^N$ with $N\ge 8$. The characters of the
bifurcation diagrams depend on different ranges of $\alp $: when
$0\le \alp \le \alp ^{**}$,
  there exists a unique solution for $(S)_{ \lam }$ with $\lam \in
(0,\lam ^*)$ and there does not
  exist any solution for $(S)_{ \lam ^*}$; when $ \alp > \alp ^{**}$,
$|u(0)|$
  oscillates around the value $\lam _*$ defined in $(\ref{1:17})$ and
there exists a unique solution for $(S)_{ \lam ^*}$.}}
\label{fig:5:1}
\end{center}
\end{figure}

\begin{figure}[htbp]
\begin{center}
{\includegraphics[width = 11cm,height=5.5cm,clip]{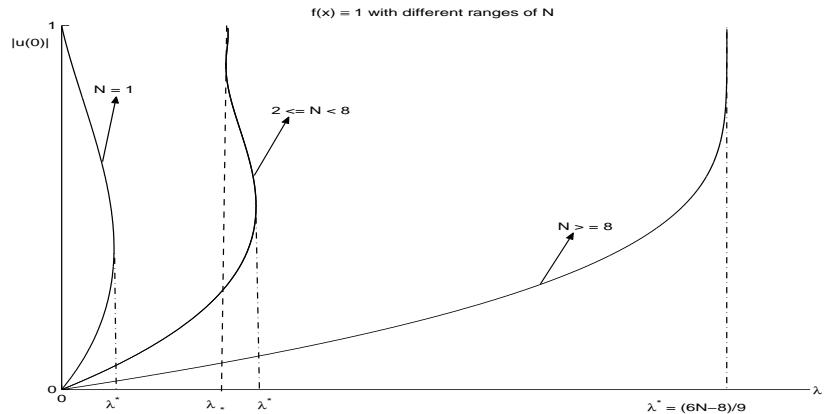}
\label{fig:fig2_2}}
  \caption{{\em   Plots of $|u(0)|$ versus $\lam$ for the constant
permittivity profile $f(x)\equiv 1$ defined in the unit ball
$B_1(0)\subset \R^N$ with different ranges of $N$. In the case of
$N\ge 8$, we have $\lam ^*=(6N-8)/9$.} } \label{fig:5:1}
\end{center}
\end{figure}

\begin{table}[htb]
\begin{center}\mbox{Number of solutions for the case}\ $f(x)=|x|^{\alp
}$\
\mbox{and}\ $N=1$\\
\begin{tabular}{  |c|c|c|c|c| }
\hline    & $\lam =\lam _1(<\lam ^*)$ & $\lam <\lam ^*$ & $\lam
=\lam ^*$
\\
\hline    $0\le \alp \le 1$  &  -----  & $2$& \mbox{1}  \\
\hline   \mbox{$1<\alp \le \alp ^{*}$}   &  $1$
& \mbox{$\ge 1$}  &\mbox{1}   \\
\hline  \mbox{$\alp > \alp ^{*}$}   &  \mbox{$\infty $}
& \mbox{$\ge 1$}  &\mbox{1}   \\
\hline \end{tabular}
\end{center}
\caption{{\em Number of solutions to $(1.1)$ which is defined in a
unit ball $B_1(0)$ with $N=1$, where $\lam _1=\frac{(2+\alp)(\alp
-1)}{9}<\lam ^*$, $ \alp ^*=-\frac{1}{2}+\frac{1}{2}\sqrt{27/2} $
and $f(x)=|x|^{\alp }$ is chosen to be power-law permittivity
profile. }} \label{tab1}
\end{table}

\vspace {.1cm}

\begin{table}[htb]
\begin{center}\mbox{Number of solutions for the case}\ $f(x)=|x|^{\alp
}$\
\mbox{and}\ $2\le N\le 7$\\
\begin{tabular}{|c|c|c|c|c|}
\hline  & $\lam =\lam _*(<\lam ^*)$ & $\lam <\lam ^*$ & $\lam =\lam
^*$
\\
\hline \mbox{$\alp \ge 0$}   &  \mbox{$\infty $}
& \mbox{$\ge 1$}  &\mbox{1}  \\
\hline \end{tabular}
\end{center}
\caption{{\em Number of solutions to $(1.1)$ which is defined in a
unit ball $B_1(0)$ with $2\le N\le 7$, where $\lam _*
=\frac{(2+\alp)(3N+\alp -4)}{9}<\lam ^*$ and $f(x)=|x|^{\alp }$ is
chosen to be power-law permittivity profile. }} \label{tab2}
\end{table} \vspace {.1cm}

\begin{table}[htb]
\begin{center}\mbox{Number of solutions for the case}\ $f(x)=|x|^{\alp
}$\
\mbox{and}\ $N\ge 8$\\
\begin{tabular}{|c|c|c|c|c|}
\hline  & $\lam =\lam _*(<\lam ^*)$ & $\lam <\lam ^*$ & $\lam =\lam
_*(=\lam ^*)$&$\lam =\lam ^*$   \\
\hline
   \mbox{$0\le \alp \le \alp ^{**}$}           &  -----  & \mbox{1}&$0$&
\mbox{0}  \\
\hline \mbox{$\alp >\alp ^{**}$}   &  \mbox{$\infty $} & \mbox{$\ge
1$} &
-----  &\mbox{1} \\
\hline \end{tabular}
\end{center}
\caption{{\em Number of solutions to $(1.1)$ which is defined in a
unit ball $B_1(0)$ with $N\ge 8$, where $\lam _*
=\frac{(2+\alp)(3N+\alp -4)}{9} $ ($=\lam ^*$ for $0\le \alp \le
\alp ^{**}$), $\alp ^{**}=\frac{4-6N+3\sqrt{6}(N-2)}{4}$ and
$f(x)=|x|^{\alp }$ is chosen to be power-law permittivity profile.}}
\label{tab3}
\end{table}

\begin{table}[htb]
\begin{center}\mbox{Number of solutions for the case}\ $f(x)\equiv 1$ \\
\begin{tabular}{|c|c|c|c|c|}
\hline  & $\lam = \lam _*(<\lam ^*)$ & $\lam <\lam ^*$ &$\lam = \lam
_*(=\lam ^*)$& $\lam =\lam ^*$   \\
\hline
   \mbox{$N=1$}         &   ----- & \mbox{   2}& -----& \mbox{ 1 }  \\
\hline \mbox{$2\le N\le 7$}   &  \mbox{$\infty $} & \mbox{$\ge 1$}
&-----
&\mbox{ 1 } \\
\hline \mbox{$N\ge 8$}           & ----- & \mbox{1 }&$0$ & \mbox{ 0
}   \\
\hline \end{tabular}
\end{center}
\caption{{\em Number of solutions to $(1.1)$ which is defined in a
unit ball $B_1(0)$, where $\lam _*=\frac{6N-8}{9}$ and $f(x)\equiv
1$ is chosen to be constant permittivity profile. We note that
$\lam_*=\lam ^*=(6N-8)/9$ for $N\ge 8$. } } \label{tab4}
\end{table}
In this case, we have $\bigtriangleup < 0$ and
$$
V\sim \Big(\frac{9}{(2+\alp)(3N+\alp -4)}\Big)^{\frac{1}{3}}+
C_1e^{-\frac{3N+2\alp -2}{6}\eta
}cos\big(\frac{\sqrt{-\bigtriangleup}}{6}\eta +C_2\big) +\cdots
\,,\quad \mbox{as} \quad \eta \to +\infty \,.
$$

Further, we conclude that
\begin{equation}
w\sim P^{\frac{2+\alp }{3}}\Big(\frac{9}{(2+\alp)(3N+\alp
-4)}\Big)^{\frac{1}{3}}+
C_1P^{-\frac{N-2}{2}}cos\big(\frac{\sqrt{-\bigtriangleup}}{6}logP
+C_2\big) +\cdots \,,\quad \mbox{as} \quad P \to +\infty \,.
\label{1:41}
\end{equation}
And we also obtain from $\lam =\gam ^{2+\alp }/w^3(\gam )$ that
$$\lam \sim \lam _*=\frac{(2+\alp)(3N+\alp -4)}{9}\quad
\mbox{as}\quad \gam \to \infty \,.$$ When $N$ and $\alp $ separately
satisfy $(\ref{1:18a})$, $(\ref{1:18b})$ and $(\ref{1:18c})$, the
typical diagrams are computed in Fig.~2, Fig.~3(a) and Fig.~3(b),
respectively.

Combining the numerical results of Fig.~2 and Fig.~3, we plot the
bifurcation diagram of the constant permittivity profile defined in
the unit ball with different ranges of $N$. The result of such a
computation is shown in Fig.~4, from which one can observe the
uniqueness, (infinite) multiplicity of the solutions for $(S)_{\lam
}$ with $\lam \in (0, \lam ^*)$ and different ranges of $N$: $N=1$,
$2\le N\le 7$ and $N\ge 8$. We note that $\lam ^*=(6N-8)/9$ when
$N\ge 8$.

Applying above numerical results, in Table 3$\sim$5 we give the
number of solutions for $(\ref{1:8})$ depending on $N$ and $\alp$.
In Table 6, we give the number of solutions for $(\ref{1:8})$ with
constant profile $f(x)\equiv 1$, which shows that $N=8$ is the
critical dimension for $(\ref{1:8})$.

\begin{rem} {\rm Under the assumptions of
$(\ref{1:18})$, since $w(P)$ in $(\ref{1:41})$ is oscillatory for $P
>>1$, we expect that $|u(0)|$ oscillates around the value $\lam
_*=\frac{(2+\alp)(3N+\alp -4)}{9}$ as $P \to \infty$. In particular,
this implies that for this case, $(S)_{ \lam}$ has infinitely
multiple solutions for $(S)_{ \lam _*}$. We compute the numerical
results of this case in Fig.~2 and Fig.~3, from which we can observe
the uniqueness, (infinite) multiplicity of the solutions for $(S)_{
\lam}$: when $N$ and $\alp $ satisfy any one of the three cases in
$(\ref{1:18})$, then there exists a series of $\{\lam _i\}$
satisfying
$$\begin{array}{lll}
\lam _0=0\,,\quad \lam _{2k}\nearrow \lam _* \quad \mbox{as}\quad
k\to \infty \,;\\[3mm]
\lam _1=\lam ^*\,,\quad \lam _{2k-1}\searrow \lam _*\quad
\mbox{as}\quad k\to \infty \,
\end{array}
$$
such that: there exist exactly $2k+1$ solutions for $(S)_{ \lam }$
with $\lam \in (\lam _{2k},\lam _{2k+2})$; and there exist exactly
$2k$ solutions for $(S)_{ \lam }$ with $\lam \in (\lam _{2k+1},\lam
_{2k-1})$; further, there exist infinitely multiple solutions for
$(S)_{ \lam _*}$. Therefore, it is reasonable to believe that the
multiplicity of solutions for the general $(S)_{\lam }$ greatly
depends on the permittivity profile $f(x)$, the dimension $N$ and
the value of $\lam $.}
\end{rem}

\begin{rem} {\rm Our results show that for
$f(x)=|x|^{\alp }$ with $N\ge 8$ and $0\le \alp \le \alp ^{**}$,
then there does not exist any classic solution for $(S)_{ \lam ^*}$,
where $\lam ^*=(2+\alp )(3N+\alp -4)/9$; but for other cases of $N$
and $\alp$, there exists a unique solution for $(S)_{ \lam ^*}$.
Therefore, for $N\ge 8$ it seems from these results that whether
there exist solutions for $(S)_{ \lam ^*}$ depends on the varying
permittivity profile $f(x)$. However, we conjecture that for $N\ge
8$ there is no solution for $(S)_{\lam ^*}$ if the permittivity
profile $f(x)>0$ on $\Omega$.}
\end{rem}

\vspace {.95cm}

\noindent {\bf Acknowledgements:} We are grateful to Michael.~J.
Ward for introducing us to the PDE models for  electrostatic MEMS
devices and for several valuable discussions concerning this paper.
We are also thankful to Louis Nirenberg for leading us to the
pioneering work of Joseph and Lundgren and the related papers of
Crandall-Rabinowitz. Special thanks also go to PierPaolo Esposito
for his thorough reading of the manuscript that led to many
improvements.


\begin{thebibliography}{GNN}




\bibitem{A} U.~Ascher, R.~Christiansen and R.~Russell, {\em Collocation
Software for Boundary Value ODE's}, Math. Comp., {\bf 33}, (1979),
pp.~659-679.



%
%
%





\bibitem {B}Bandle, C., {\em Isoperimetric Inequalities and Applications},
  In Monographs and Studies in Mathematics, Boston,
Mass.-London, Pitman (1980).



\bibitem{B1} R.~E. Bank, {\em PLTMG: A Software Package for Solving
Elliptic Partial Differential Equations}, User's guide 8.0,
Software, Environments, and Tools, SIAM, Philadelphia, PA, (1998),
xi+110 pages.







\bibitem {BGP} D. Bernstein, P. Guidotti and J.~A. Pelesko,
{\em Analytic and numerical analysis of electrostatically actuated
MEMS devices}, Proc. of Modeling and Simulation of Microsystems,
{\bf 2000}, (2000), pp.~489--492.

%

\bibitem {CM}X.~Cabre and Y.~Martel, {\em Weak eigenfunctions for the
linearization of extremal elliptic problems}, J. Funct. Anal., 156
(1998), pp.~30--56.

\bibitem {CR1} M.~G.~Crandall and P.~H.~Rabinowitz, {\em Bifurcation,
perturbation of simple eigenvalues and linearized stability}, Arch.
Rational Mech. Anal., 52 (1973), pp.~161--180.

\bibitem {CR} M.~G.~Crandall and P.~H.~Rabinowitz, {\em Some continuation and
variational methods for positive solutions of nonlinear elliptic
eigenvalue problems},  Arch. Ration. Mech. Anal., 58 (1975),
pp.~207--218.


\bibitem {E} L.~Evans, {\em Partial differential equations}, Graduate Studies
in Mathematics, 19. AMS, Providence, RI, 1998.


\bibitem {FK} S. Filippas and R.~V. Kohn, {\em Refined Asymptotics for
the Blow Up of $u_t-\Delta u=u^p$}, Comm. Pure Appl. Math., {\bf
45}, No. 7, (1992), pp.~821--869.



\bibitem {FMP} G. Flores, G.~A. Mercado and J.~A. Pelesko, {\em Dynamics
and Touchdown in Electrostatic MEMS}, Proceedings of ICMENS 2003,
(2003), pp.~182--187.


\bibitem {GG2}N. Ghoussoub and  {Y. Guo}, {\em On the partial
differential equations of electrostatic MEMS devices II: dynamic
case}, in preparation.



\bibitem {GG3}N. Ghoussoub and  {Y. Guo}, {\em On the partial
differential equations of electrostatic MEMS devices III: refined
touchdown behavior}, in preparation.




\bibitem {GNN} B. Gidas, W.~M. Ni and L. Nirenberg, {\em Symmetry and
related properties via the maximum principle}, Comm. Math. Phys.,
{\bf 68}, (1979), pp.~209--243.





\bibitem {GPW} Y. Guo, Z. Pan and Michael J. Ward,
{\em Touchdown and pull-in voltage behavior of a MEMS device with
varying dielectric properties}, accepted by SIAM, J. Appl. Math.



\bibitem {HW} A. Haraux and F.~B. Weissler, {\em Non-uniqueness for a
semilinear initial value problem}, Indiana Univ. Math. J., {\bf 31},
No. 2, (1982), pp.~167--189.



\bibitem {JL} D.~D.~Joseph and T.~S.~Lundgren, {\em Quasilinear Dirichlet
problems driven by positive sources}, Arch. Ration. Mech. Anal.,
{\bf  49} (1973), pp.~241--268.



\bibitem {KK}J.~P.~Keener and H.~B.~Keller, {\em Positive solutions
of convex nonlinear eigenvalue problems}, J. Diff. Eqns, 16, (1974),
pp.~103--125.






\bibitem {NN} H.~C. Nathanson, W.~E. Newell and R.~A. Wickstrom, J.~R.
Davis, {\em The Resonant Gate Transistor}, IEEE Trans. on Elect.
Devices, {\bf 14}, (1967), pp.~117--133.



\bibitem {P1} J.~A. Pelesko, {\em Mathematical Modeling of
Electrostatic MEMS with Tailored Dielectric Properties}, SIAM J.
Appl. Math., {\bf 62}, No. 3, (2002), pp.~888--908.







\bibitem {PB} J.~A. Pelesko and D.~H. Bernstein, {\em Modeling MEMS and
NEMS}, Chapman Hall and CRC Press, (2002).



\bibitem {PT} J.~A. Pelesko and A.~A. Triolo, {\em Nonlocal Problems in
MEMS Device Control}, J. Eng. Math., {\bf 41}, No. 4, (2001),
pp.~345--366.



\bibitem {S} I. Stackgold, {\em Green's Functions and Boundary Value
Problems}, Wiley, New York, (1998).





\bibitem {SA} M.~T.~A. Saif, B.~E. Alaca and H. Sehitoglu,
{\em Analytical modeling elctrostatic membrane actuator micropumps},
IEEE J. of MEMS, No. 8, (1999), pp.~335--344.



\bibitem {T} G.~I. Taylor, {\em The coalescence of closely spaced
drops when they are at different electric potentials}, Proc. Roy.
Soc. A, {\bf 306}, (1968), pp.~423--434.













\end{thebibliography}
\end{document}